\documentclass[12pt]{article}
\usepackage{amsmath, amscd, amssymb, latexsym,hhline}
\usepackage[all]{xy}
\def\:{{\colon}}
\def\sgn{{\mbox{sgn}}}
\def\End{{\mbox{\rm End}}}
\def\Ext{{\mbox{\rm Ext}}}
\def\Tor{{\mbox{\rm Tor}}}

\def\Hom{{\mbox{\rm Hom}}}
\def\Aut{{\mbox{\rm Aut}}}

\def\ker{{\mbox{ker\,}}}
\def\coker{{\mbox{coker\,}}}

\def\Im{{\mbox{Im\,}}}
\def\ord{{\mbox{ord}}}

\def\BC{{\mathbb C}}
\def\BZ{{\mathbb Z}}
\def\BR{{\mathbb R}}
\def\BQ{{\mathbb Q}}

\def\infl{{infl\,}}

\def\a{{\alpha}}
\def\b{{\beta}}
\def\d{{\delta}}
\def\l{{\lambda}}
\def\g{{\gamma}}
\def\G{{\Gamma}}

\def\Hone{{D^{\w_1}(G)}}
\def\H'{{D^{\w'}(G)}}
\def\H{{D^\w(G)}}
\def\Htwo{{D^{\w_2}(G)}}
\def\Hi{{D^{\w_i}(G)}}

\def\e{{\epsilon}}
\def\s{{\sigma}}
\def\t{{\tau}}

\def\w{{\omega}}
\def\L{{\Lambda}}

\def\p{{\partial}}
\def\map{{\longrightarrow}}
\def\qed{{\ \ \ $\square$}}
\def\pf{{\it Proof\ \ }}

\def\Wedge{{\bigwedge\nolimits}}
\def\<{{\langle}}
\def\>{{\rangle}}
\def\pf{{\it Proof. }}
\def\C{{\,\subseteq\,}}

\newcommand{\M}[1]{{#1\mbox{-\bf Mod}}}
\newcommand{\DG}[1]{{D^{\w_{#1}}(G_{#1})}}

\newtheorem{thm}{Theorem}[section]
\newtheorem{cor}[thm]{Corollary}
\newtheorem{defn}[thm]{Definition}
\newtheorem{prop}[thm]{Proposition}
\newtheorem{remark}[thm]{Remark}

\newtheorem{lem}[thm]{Lemma}

\newtheorem{example}[thm]{Example}

\begin{document}
\title{Group Cohomology and Gauge Equivalence of some Twisted Quantum Doubles}
\author{Geoffrey Mason\footnote{Research supported by the National Science Foundation and the
Regents of the University of California.} and Siu-Hung Ng \\ Mathematics Department, University
of California\\ Santa Cruz, CA 95064}
\date{December 5, 1999}
\maketitle

\section{Introduction}
   The purpose of the present paper is to investigate relationships
that exist between certain  quasi-Hopf algebras (namely, twisted quantum doubles of a finite
abelian group) and their module categories. This is closely related to
aspects of group cohomology concerning the relationship of the bar complex to some other
complexes defined and studied by Eilenberg and MacLane \cite{Mac52}, \cite{EMac1} and
\cite{EMac2}. Although our results are purely homological and algebraic in nature, the
motivation for studying these questions derives mainly from connections with conformal field
theory (CFT) and vertex operator algebras
\cite{DW}, \cite{DPR90}, \cite{MS}, \cite{DL}. We will thus review background motivation (although
it is hardly necessary for an understanding of our main results) and then give a more
detailed overview of the paper.

   It is not necessary for the reader to appreciate what a vertex operator algebra
(VOA) is (cf. \cite{Borcherds}, \cite{FLM}, \cite{FHL}), but only to understand this: VOAs have a representation
theory, hence a (linear) module category, and it is important to understand the
nature of this category. It is expected that the module category $\M{V}$ for a VOA $V$ has the
structure of a braided monoidal category \cite{JS93}. Roughly, this means that there is a notion of
tensor product of $V$-modules, and that the tensor product of three modules is associative by
an isomorphism that depends on the particular modules. Moreover, tensor product of modules is
commutative, again via an isomorphism that depends on the particular modules. Particularly
nice VOAs called {\em rational} have the additional property that $\M{V}$ is a semi-simple
category, so that $V$ has only a finite number of isomorphism classes of simple modules, and
every object in $\M{V}$ is a direct sum of simple objects. The simplest VOAs of all from this
perspective are called {\em holomorphic}: they are rational and have but one simple module,
namely the adjoint module $V$ itself. Thus if $V$ is a holomorphic VOA, then $\M{V}$ is
equivalent (as a semi-simple braided monoidal category) to the category of vector spaces over
$\BC$. Examples of holomorphic VOAs include the famous Moonshine
Module  and VOAs attached to positive-definite, even, self-dual lattices \cite{FLM}.

   It is a very difficult problem to establish the nature of $\M{V}$, even for
relatively well-understood VOAs. However, if $\M{V}$ is indeed a braided monoidal category, then
reconstruction theory (cf. \cite{haring}, \cite{Maji93} ) suggests that there might be a
quasi-Hopf algebra
$Q$ with the property that $\M{Q}$ is equivalent ( as braided monoidal category) to
$\M{V}$. (Recall that the module category of a quasi-Hopf algebra automatically carries a
braided monoidal structure, cf. \cite{Kassel}, \cite{Drin90}.) In this context, Dijkgraaf,
Pasquier and Roche
\cite{DPR90} put forward a remarkable proposal which essentially identifies such a $Q$ in the case
that $V$ is a so-called {\em holomorphic orbifold}. This means that there is a holomorphic VOA
$W$ and a finite group of automorphisms $G$ of $W$ such that $V = W^G$, i.e., $V$ is the subVOA
of
$G$-fixed-points. The DPR ansatz is that one can take $Q$ to be a twisted quantum double $\H$
of $G$. This is a certain quasi-Hopf algebra, first constructed in \cite{DPR90}, which is a version of
the Drinfeld double of the group algebra $\BC[G]$ \cite{DrinICM}, but which is twisted by a
certain 3-cocycle $\w$ in $Z^3(G,\BC^*)$. For more details see \cite{DPR90}, \cite{Mas}, and
section \ref{s1} below.

   One of the problems of vertex operator algebras is that of classification. There
are far too many VOAs to make the problem of classification up to isomorphism at all
practical, at least for the foreseeable future. A meaningful alternative is to classify VOAs
up to equivalence of their module categories. Thus we say that two VOAs $V$ and $W$ are
{\em physically equivalent}  if $\M{V}$ and $\M{W}$ are equivalent as linear, braided monoidal
categories. For example, all holomorphic VOAs are equivalent in this sense because they have
identical module categories. In the present paper we are interested in the case of
holomorphic orbifolds: given two holomorphic VOAs $W_1$ and $W_2$, and finite automorphism
groups $G_1$ and $G_2$ of $W_1$ and $W_2$ respectively, when are the corresponding holomorphic
orbifolds physically equivalent? Because of the DPR ansatz, one expects that this is
equivalent to the following purely algebraic question: \\

\hspace{0.5cm} (*)\hspace{0.5cm}
\begin{minipage}[c]{5in}
given finite groups $G_1$ and $G_2$ and normalized 3-cocycles $\w_1 \in Z^3(G_1,\BC^*)$ and
              $\w_2 \in Z^3(G_2,\BC^*)$, when are the module categories $\M{D^{\w_1}(G_1)}$
and
              $\M{D^{\w_2}(G_2)}$ equivalent as (linear) braided monoidal categories?\\
\end{minipage}

   As far as we know, the question in this generality has not been
discussed. In the present paper we
will limit ourselves to the case in which the relevant groups and
algebras are commutative. Note that it is well-known \cite{Kassel} that a
sufficient condition for the equivalence of the categories in (*) is
that the corresponding quasi-bialgebras are gauge equivalent. We will
prove (Theorem \ref{t8.5}) that this is also a necessary condition, at least
when the twisted quantum doubles are commutative. In practice, however,
this is not very helpful in trying to decide the answer to (*) in a
given case. The main results and methods of this paper show essentially
that there is a method that allows one to answer (*) in any given
case. We also establish results that give a complete answer for various
classes of groups.

   Let us begin by considering a finite abelian group $G$. With each normalized 3-cocycle $\w$
in $Z^3(G,\BC^*)$ we consider the twisted quantum double $\H$ together with its group of
group-like elements $\G^\w(G)$. Now $\G^\w(G)$ is an abelian group, and we will be
particularly interested in the case in which $\G^\w(G)$ spans $\H$: this happens precisely
when $\H$ is isomorphic as a bialgebra to the group algebra $\BC[\G^\w(G)]$, in which case
$\G^\w(G)$ is a central extension of $G$ by its group of characters $\widehat{G}$.
\begin{equation}\label{eqi1}
                        1\map \widehat{G}\map \Gamma^\w(G) \map G
\map 1
\end{equation}

However,
there is an essential difference between the two bialgebras which generally means that they
are not isomorphic as quasi-bialgebras: the associator of $\H$ carries information about the
cocycle $\w$, whereas the associator of the group algebra is trivial and carries no
cohomological information.

   We call the normalized 3-cocycle $\w$ {\em abelian} in case $\G^\w(G)$ spans $\H$, and denote
by $Z^3(G,\BC^*)_{ab}$, resp $H^3(G,\BC^*)_{ab}$ the 3-cocycles, resp. 3-cohomology classes
which are abelian in this sense. They are groups, and via (\ref{eqi1}) we get a group
homomorphism
\begin{equation}\label{eqi2}
 \Lambda : H^3(G,\BC^*)_{ab}
\map H^2_{ab}(G,\widehat{G})
\end{equation}
 where the subscript {\em ab} on the second cohomology group refers to {\em abelian}
extensions of
$G$ by $\widehat{G}$.

   It transpires that many of our results hinge on a close analysis of this map. Though
we have described $\L$ in terms of twisted quantum doubles, it turns out to have a purely
homological interpretation. This arises as follows: Eilenberg and MacLane in \cite{Mac52}, \cite{EMac1} and
\cite{EMac2} described a sequence of complexes $A_i (G)$ associated to an abelian group $G$ with
a view to computing the cohomology of the Eilenberg-MacLane spaces $K(G,n)$. We are only
interested in the complex $A_0 (G)$, which is just the bar complex, and the second complex
$A_1 (G)$ which contains $A_0 (G)$ as a subcomplex. This containment yields a canonical short
exact sequence of complexes
$$
0 \map \Hom(B(G), \BC^*) \map \Hom (A_1(G),\BC^*) \map \Hom(A_0(G), \BC^*) \map 0
$$
($B(G)$ is the quotient complex  $A_1(G) / A_0(G)$), and thereby a long
exact sequence in cohomology. It turns out that the map $\L$ in (\ref{eqi2}) is the
restriction to
$H^3(G,\BC^*)_{ab}$ of the connecting homomorphism which maps $H^3(G,\BC^*)$ to
$H^4(B(G),\BC^*)$. There is naturally a close connection between $\L$ and the cohomology of
the complex $A_1(G)$, which is sometimes called abelian cohomology.  We observe that this latter
cohomology already plays a role in
\cite{MS} and in
\cite{JS93}, and is quite essential in the work of Dong and Lepowsky \cite{DL} on the
generalized Jacobi identity satisfied by the vertex operators associated to lattice elements in
a lattice VOA. But the context of these works is quite different to our own.

 Our analysis of the long exact sequence also reveals a fundamental
difference between groups of odd order and groups of even order : we
will see that $\L$ is injective if, and only if, $G$ has odd order, and that in general $\ker\L$
is isomorphic to the group\footnote{For a finite abelian group $G$, $\Omega_p(G)$ is the
subgroup of elements of order dividing $p$.}
$\Omega_2(G)$ of elements in
$G$ of order at most 2. We also give a precise description of
$\Im
\L$ : in the easier case of groups of odd order it is precisely the group of invariants of the
duality map $\e$ : $H^2_{ab}(G, \widehat{G}) \map H^2_{ab}(G, \widehat{G})$
which maps a
short exact sequence of abelian groups to the dual short exact sequence
obtained by applying the functor $\Hom( ?, \BC^*)$. In case $G$ is an abelian
2-group, it will turn out from the structure of $\Im \L$ that the
extension $\G^\w(G)$ is necessarily a product of an {\em even}
number of cyclic factors. This and other results in the paper
have consequences for the theory of VOAs and CFT that we hope to pursue
elsewhere.

   The main consequence of our analysis of the long exact sequence is
that we are able to characterize the pair $(G,\w)$ for $\w$ abelian in terms
that are more amenable to calculation. Namely, we show that it is
equivalent to the existence of a certain non-degenerate quadratic form
on the abelian group $\G^\w(G)$, call it $q$, with the property that $q$
has a {\em metabolizer}  $\widehat{G}$. In general we say that a
non-degenerate quadratic space $(E,q)$ for a finite abelian group $E$ has a
metabolizer $M$, if the restriction of $q$ to $M$ vanishes identically and if
the order of $E$ is the square of the order of $M$.
So here we are following the definition of metabolizer in \cite{T} rather than
the classical definition \cite{Knebusch}, \cite{HM},
where $M$ is also required to be a direct summand of $E$. In the case we
are interested in, this only happens if the extension (\ref{eqi1}) splits,
and we say that $M$ is a {\em split metabolizer} in this case. For us, metabolizers are
not necessarily split and the most interesting metabolizers are not. In this context we may speak
of metabolic triples $(E, q, G)$, and we can try to construct suitable Witt-type
groups. We  essentially carry this out when $E$ is a homogeneous (or
homocyclic) abelian $p$-group, however we do not use the language of Witt groups.
We will see
(Theorem
\ref{t8.4}) that the equivalence (*) for abelian groups $G_1$, $G_2$ and abelian 3-cocycles 
$\w_1, \w_2$ is the same as the equivalence of the associated quadratic spaces
$(\G^{\w_1}(G_1), q_1)$ and $(\G^{\w_2}(G_2), q_2)$.

   We point out that quadratic forms play a prominent role in several
of the papers we have already quoted concerning abelian cohomology and/or
VOAs (\cite{DL}, \cite{JS93}, \cite{MS}). But in \cite{DL} and \cite{MS}, for example, conformal
field theories more general than holomorphic orbifolds are considered, so that the existence of a
metabolizer, so crucial for us, plays no role in these works.

   Once the connection with quadratic forms with a metabolizer has been
forged, a number of further results flow. Consider a finite abelian
group $G$ together with an abelian 3-cocycle $\w$. By analogy with the
standard decomposition of $G$ into a direct product of cyclic groups, one
may ask whether $D^\w(G)$ is gauge equivalent to a tensor product
\begin{equation}\label{eqi3}
D^\w(G) \sim D^{\w_1}(C_1)  \otimes \cdots \otimes
D^{\w_r} (C_r)
\end{equation}
with each $C_i$ a cyclic group. One may assume that $G$ is a $p$-group for
some prime $p$, in which case the question has an affirmative answer in
the following cases: (a)  $p$ satisfies the congruence  $p \equiv 1 \pmod{4}$;
(b)  $G$ is a homogeneous $p$-group and $p$ is odd. However, the result is
generally false if $p = 2$, and $p \equiv 3 \pmod{4}$, and $G$ is not
homogeneous. Furthermore, if there is an equivalence of the form (\ref{eqi3})
then it is generally not unique, and the direct product   $C_1 \times \cdots \times
C_r$  of the groups $C_i$  is not necessarily isomorphic to $G$.  It may
seem curious that one has to impose congruence conditions on $p$ in
order to establish (\ref{eqi3}) for all abelian $p$-groups. The reason is that
via the quadratic forms we can use the technique of Gauss sums, and odd primes
which are sums of two squares are better behaved.  The proof of (\ref{eqi3}) for homogeneous
p-groups is established using Witt-group type techniques. The paper of C.T.C. Wall \cite{Wa63},
which classifies the quadratic forms on finite abelian groups (at least
if they are of odd order) proves to be invaluable in many of our
calculations.

By a {\em rational lattice}  we mean a free abelian group $M$ of finite
rank together with a non-degenerate bilinear form  $\<\bullet , \bullet \> :  M \times M
\map \BQ$. We do not necessarily  assume that  $\<\bullet , \bullet \>$  is
positive-definite. The lattice is called {\em integral}  in case it takes
values in $\BZ$, and {\em even} if, in addition, one has  $\< x , x \>  \in 2\BZ$. Such
lattices play a prominent role in \cite{DL} because they naturally give
rise to vertex algebras, indeed to examples of what Dong-Lepowsky
call ``abelian intertwining algebras.'' Via the connection with
quadratic forms, we are able to relate the elements of $H^3(G,\BC^*)_{ab}$
to certain pairs of lattices  $M \C L$, essentially because such a pair
(together with the bilinear form  $\<\bullet , \bullet \>$  ) yields a quadratic form
on $L / M$. Moreover, a theorem of Wall \cite{Wa63} allows us to see that the
existence of a metabolizer implies that the lattice $M$ is necessarily
self-dual (i.e., unimodular). As we will show elsewhere, this has
consequences for the theory of lattice VOAs. It also permits us to
see that the following pieces of data are essentially equivalent\footnote{a vague phrase which
will be clarified in section \ref{s10}}  for a finite abelian group $G$:
\begin{center}
\begin{enumerate}
            \item[(a)] $[\w] \in H^3(G, \BC^*)_{ab}$.
            \item[(b)] a metabolic triple $(E, q, \widehat{G})$.
            \item[(c)] a pair of rational lattices  $M \C L$  with $M$ even and
self-dual and  $L / M  =  G$.
\end{enumerate}
\end{center}
    We will also prove (Theorem \ref{t11.6}) that if $G$ is a homogeneous abelian $p$-group of
odd order then two twisted quantum doubles $D^\w(G)$ and $D^{\w'}(G)$ for abelian cocycles $\w$
and
${\w'}$ are gauge equivalent if, and only if, the corresponding cohomology classes $[\w]$ and
$[{\w'}]$ are equivalent under the action of the automorphism group $\Aut(G)$. In some sense
this result is the quintessential goal that one seeks for any abelian
$p$-group $G$; but the result is false in general due to ``hidden'' gauge
equivalences which are hard to enumerate, but which can be detected
via the metabolic triples.

    Finally, the connection with quadratic forms allows us to
construct various kinds of dualities and symmetries (other than
categorical equivalence) between the module categories in (*). We
give various examples in section \ref{s13}, including the following: for
positive integers $n$, $k$ and odd prime $p$, let $G_1$ and $G_2$ denote
the homogeneous groups $\left(\BZ_{p^n}\right)^k$ and $\left(\BZ_{p^k}\right)^n$ respectively.
Then  there are precisely ${{n + k} \choose k}$
equivalence classes of monoidal categories (i.e., tensor categories
excluding braiding) both of the form $\M{D^{\w_1}(G_1)}$ and $\M{D^{\w_2}(G_2)}$, for some
$\w_1$
or
$\w_2$ respectively. Moreover there is a  canonical bijection (duality) between these two sets of
equivalence  classes of tensor categories. Each such tensor category can be
labeled by a partition, and at the level of partitions the duality
is simply that which maps a partition to its dual (or conjugate)
partition. Moreover, duality induces a bijection between the braided
monoidal categories associated to each group.

    The paper is organized as follows: we give some background in
section \ref{s1}, including the observation that for a finite abelian group
and normalized 3-cocycle $\w \in Z^3(G, \BC^*)$, the object $\H$--when equipped
with the trivial associator--is a semi-simple self-dual Hopf algebra
which is generally neither commutative nor cocommutative. We study
the group-like elements of $\H$ in section \ref{s2}, giving rise to the
central extension (\ref{eqi1}). After a short section \ref{s3} concerned with
tensor products of twisted quantum doubles, we take up in sections \ref{s4} - \ref{s7} the
structure of the long exact sequence alluded to above. Much
of the resulting homological information is encoded in a 7 term exact
sequence which we record in Remark \ref{r6.6}. In section \ref{s8} we consider the
category $\M{\H}$  as a monoidal category i.e., forgetting the
braiding. For example, we show (Lemma \ref{r7.1} and Theorem \ref{t7.3}) that for
finite abelian groups $G$ of odd order and normalized 3-cocycles $\w_1$ and $\w_2$ in
$Z^3(G, \BC^*)_{ab}$, the monoidal categories $\M{D^{\w_1}(G)}$ and $\M{D^{\w_2}(G)}$
are tensor equivalent if, and only if, the corresponding groups
of fusion rules $\G^{\w_1}$ and $\G^{\w_2}$ are isomorphic, and that
in this case each $D^{\w_i}(G)$ is gauge equivalent to the group algebra
$\BC[\G^{\w_i}]$. This answers in the affirmative--at least for groups
of odd order--a question posed in \cite{DPR90}, namely can $\H$ be
obtained by twisting a Hopf algebra? We show by example that the
answer is ``no'' in general. In section \ref{s9} we study the connections
among gauge equivalence, braided monoidal categories and quadratic
forms and prove the results already mentioned. Section \ref{s10} explains
the connections with lattices, while section \ref{s11} presents
some results involving Gauss sums and applies them to the proof of
the gauge equivalence
(\ref{eqi3}) in the case $p \equiv 1 \pmod{4}$. In section \ref{s12}
we consider metabolic triples $(E, q, G)$ when $G$ is a homogeneous
$p$-group and--in all but name--construct the corresponding Witt-type
group. This is used (Theorem \ref{t11.6}) to establish both (\ref{eqi3}) in the
case that $G$ is a homogeneous abelian $p$-group, and other results concerning homogeneous
$p$-groups already discussed. Finally, we discuss symmetries and
dualities in section \ref{s13}.

    The authors thank Chongying Dong for useful discussions.
\section{Twisted Quantum Double Of A Finite Group}\label{s1}
Let $G$ be a finite group and $\w: G \times G \times G \map \BC^*$ be a normalized
3-cocycle\footnote{All cocycles will take values in a trivial $G$-module.}. For any $x,y,g \in G$, define
\begin{eqnarray}
\theta_g(x,y) &=&\frac{\w(g,x,y)\w(x,y,(x y)^{-1}g x y)}{\w(x,x^{-1}g x,y)}\label{eq0.001}\\
\g_g(x,y) &=& \frac{\w(x,y,g)\w(g, g^{-1}x g, g^{-1}yg)}{\w(x,g, g^{-1}y g)}\label{eq0.002}
\end{eqnarray}
The {\em twisted quantum double} $D^\w(G)$ of $G$ with respect to $\w$ is the
quasi-triangular quasi-Hopf algebra with underlying vector space $(\BC G)^* \otimes \BC G$
and multiplication,
\begin{equation}\label{multiplication}
(e(g) \otimes x)(e(h) \otimes y) =\theta_g(x,y) e(g)e(xhx^{-1})\otimes x y
\end{equation}
\begin{equation}\label{comultiplication}
\Delta(e(g)\otimes x)  = \sum_{hk=g} \g_x(h,k) e(h)\otimes x \otimes e(k) \otimes x
\end{equation}
\begin{equation}\label{eq0.01}
\Phi = \sum_{g,h,k \in G} \w(g,h,k)^{-1} e(g) \otimes 1 \otimes e(h) \otimes 1 \otimes
e(k) \otimes 1
\end{equation}
where $\{e(g)|g \in G\}$ is the dual basis of the canonical basis of $\BC G$ (cf.
\cite{DPR90}). The counit and antipode are given by
$$
\e(e(g)\otimes x) = \delta_{g,1} \quad \mbox{and} \quad
S(e(g)\otimes x) =
\theta_{g^{-1}}(x,x^{-1})^{-1}\g_x(g,g^{-1})^{-1}e(x^{-1}g^{-1}x)\otimes x^{-1}
$$
where $\delta_{g,1}$ is the Kronecker delta. The corresponding elements $\a$ and $\b$
are $1_{D^{\w}(G)}$ and $\sum\limits_{g \in G} \w(g,g^{-1},g)e(g)\otimes 1$ respectively.
For the definition and more details about quasi-Hopf algebras, see
\cite{Drin90},
\cite{Kassel} or
\cite{Chari}.
Verification of the detail involves the following identities, which result
from the 3-cocycle identity for
$\w$:
\begin{equation}\label{eq0.1}
\theta_z(a,b)\theta_z(a b,c) = \theta_{a^{-1}z a}(b,c)\theta_z(a,b c)\,,
\end{equation}
\begin{equation}\label{eq0.2}
\theta_y(a,b)\theta_z(a,b)\g_a(y,z)\g_b(a^{-1}y a,a^{-1}z a)=\theta_{yz}(a,b)\g_{ab}(y,z)\,,
\end{equation}
\begin{equation}\label{eq0.3}
\g_z(a,b)\g_z(a b,c)\w(z^{-1}a z,z^{-1}b z,z^{-1}c z)=\g_z(b,c)\g_z(a,b c)\w(a,b,c)\,,
\end{equation}
for all $a,b,c,y,z \in G$. The universal $R$-matrix is given by
\begin{equation}\label{eq0.31}
R = \sum_{g,h \in G} e(g) \otimes 1 \otimes e(h) \otimes g \,.
\end{equation}
\begin{remark}\label{r0.1}\hspace{1cm}
{\rm
\begin{enumerate}
\item[(i)] Any Hopf algebra $H$ can be viewed as a quasi-Hopf algebra with the trivial
associator $\Phi=1_H \otimes 1_H \otimes 1_H$ and $\a=\b=1_H$.
\item[(ii)] If $w$ and $w'$ are  cohomologous 3-cocycles, then $D^\w(G)$ and $D^{\w'}(G)$
are {\em gauge equivalent  as quasi-triangular quasi-bialgebras} or  simply {\em
gauge equivalent} (cf.
\cite{Kassel}) with the algebra  isomorphism $\Theta\: D^\w(G) \map D^{\w'}(G)$ defined by
$$
\Theta(e(g) \otimes x) = \frac{b(g,x)}{b(x,x^{-1}gx)} e(g) \otimes x
$$
and the {\em gauge transform}
$F \in D^{\w'}(G) \otimes D^{\w'}(G)$ given by
$$
F=\sum_{g,h} b(g,h)^{-1} e(g) \otimes 1 \otimes e(h) \otimes 1
$$
where $\w'=  \w \delta b$ (cf. \cite{DPR90}).
\item[(iii)] Let $\s$ be a group automorphism of $G$. Then $\widehat{\s}\w$ defined by
$$ \widehat{\s}\w(g,h,k) = \w(\s^{-1}(g),\s^{-1}(h),\s^{-1}(k)) $$ is a
3-cocycle. Moreover, $\s$ induces an isomorphism $\tilde{\s}$ of
quasi-triangular quasi-Hopf algebras from $D^\w(G)$ to
$D^{\s\w}(G)$, namely $$ \tilde{\s}(e(g) \otimes x) = e(\s(g))
\otimes \s(x)\,.
$$
\end{enumerate}
}
\end{remark}

For $g \in G$, denote by $C_G(g)$ the centralizer of $g$ in $G$. One can easily see that
$\theta_g$ is a 2-cocycle of $C_G(g)$ with coefficient in $\BC^*$. Moreover, the map
$D_g \:\w \mapsto \theta_g$ induces a group homomorphism from $H^3(G,\BC^*)$ to
$H^2(C_G(g),\BC^*)$. We will denote by
$H^3(G,\BC^*)_{ab}$ the group $\bigcap\limits_{g \in G} \ker D_g$ and
$Z^3(G,\BC^*)_{ab}$ the subgroup of normalized 3-cocycles whose
cohomology classes are in $H^3(G,\BC^*)_{ab}$.\\

When $G$ is abelian, $\theta_g = \g_g$ and $C_G(g) =G$ for all $g
\in G$ and we will simply write $\theta_g$ as $\w_g$.
It follows from equation (\ref{eq0.2}) that we have
$$
\dfrac{\w_g (x,y)\w_h(x,y)}{\w_{gh}(x,y)} =
\dfrac{\w_{xy}(g,h)}{\w_x(g,h)\w_y(g,h)}\,.
$$
Hence the map $G \map H^2(G, \BC^*)$, $g \mapsto [\w_g]$ is a group homomorphism.
Let us denote this map by $\Omega(\w)$. Notice that $\Omega(\w) \in H^1(G,
H^2(G,\BC^*))$ and $\Omega(\w\w')= \Omega(\w)\Omega(\w')$ for any other 3-cocycle $\w'$.
Moreover, $\Omega(\w)$ is independent of the choice of representative of the cohomology
class of $\w$. Therefore, $\Omega$ induces a group homomorphism $\overline{\Omega}$ from
$H^3(G,\BC^*)$ to $H^1(G, H^2(G,\BC^*))$, and $H^3(G,\BC^*)_{a b} = \ker
\overline{\Omega}$.\\

\begin{prop}\label{p0.4}
Let $G$ be a finite abelian group and $\w$ a normalized 3-cocycle.
\begin{enumerate}
  \item[\rm (i)]When equipped with the trivial associator,
  $D^\w(G)$ is a self-dual Hopf algebra.
  \item[\rm (ii)] Let $b: G \times G \map \BC^*$ be a normalized 2-cochain and $\w' = \w \delta b$.
  Then the map $\Theta : D^\w(G) \map D^{\w'}(G)$ defined by
  $$
 \Theta(e(g) \otimes x ) =  \frac{b(g,x)}{b(x,g)} e(g) \otimes x
  $$
  is an isomorphism of Hopf algebras.
\end{enumerate}
\end{prop}
\pf (i) Since $G$ is abelian, $e(g) \otimes 1$ is in the center of
$D^\w(G)$. Hence, $\b \in \mbox{Cent}(D^\w(G))$ and
$\Phi \in \mbox{Cent}(\H^{\otimes 3})$. Therefore, $\Delta$
is coassociative and $S$ satisfies the antipode conditions.\\
To show that $D^\w(G)$ is self-dual, let us denote by $\{f_{x,g}\}_{g,x \in G}$
 the dual basis of $\{e(g) \otimes x \}_{g,x \in G}$.  Then,
\begin{eqnarray*}
f_{x,g} \ast f_{y,h} (e(k) \otimes z) & =& (f_{x,g} \otimes f_{y,h})\left(\sum_{ab=k}
\w_z(a,b)e(a)\otimes z \otimes e(b) \otimes z\right)\\
&=& \delta_{x,z} \delta_{y,z} \delta_{gh,k} \w_x(g,h) \\
&=& \delta_{x,y}\w_x(g,h) f_{x,gh}(e(k) \otimes z) \,.
\end{eqnarray*}
Therefore,
$$
f_{x,g} \ast f_{y,h} = \delta_{x,y}\w_x(g,h) f_{x,gh}\,.
$$
Similar, one can derive that
$$
\Delta(f_{x,g}) = \sum_{uv=x} \w_g(u,v) f_{u,g} \otimes f_{v,g}\,.
$$
The identity, $1_{\H^*}$, and the counit, $\epsilon$, of $D^\w(G)^*$ are given by
$$
1_{\H^*}=\sum_{x \in G} f_{x,1}, \quad \epsilon(f_{x,g})=\delta_{x,1}\,.
$$
respectively, Therefore, the linear map $\varphi: \H \map \H^*$ defined by
\begin{equation}\label{eq0.4}
 \varphi(e(x) \otimes g)=f_{x,g}
\end{equation}
is a Hopf algebra isomorphism. \\

(ii) Notice that for any $x,y, g \in G$
$$
\w'_g(x,y) = \w_g(x,y) \delta b_g(x,y)
$$
where $b_g(x) = \frac{b(x,g)}{b(g,x)}$. The result follows immediately from this
identity. \qed\\

In the sequel, we will denote by $\H_0$ the Hopf algebra structure on $\H$, and
we use $\H$ for the quasi-Hopf algebra structure with respect to the associator
$\Phi$ given by the equation (\ref{eq0.01}).
%%%%%%%%
%
%
%\input section2.tex
\section{Group-like Elements of $D^\w(G)$}\label{s2}
\begin{defn}{\rm
Let $G$ be a finite group and $\w$ a normalized 3-cocycle.
A nonzero element $u$ in $D^\w(G)$ is called {\em group-like} if $\Delta(u) = u \otimes
u$. We will denote by $\G^\w(G)$ the set of all group-like elements of $D^\w(G)$.
We  simply write $\G^\w$ for $\G^\w(G)$ when the context is clear.
}
\end{defn}
As in the case of group-like elements in a coalgebra (cf. \cite{Sw69}),
$\G^\w$ is a linearly independent set and any group-like element
$u$ is invertible in $D^\w(G)$ with inverse
$S(u)$. Since
$\Delta : D^\w(G) \map  D^\w(G) \otimes D^\w(G)$ is an algebra map, $\G^\w$ is a
subgroup of the group of units of $D^\w(G)$. Moreover,
$\G^\w$ can be characterized by the following proposition.
\begin{prop}\label{p0.1}
A nonzero element $u$ in $D^\w(G)$ is a group-like element if, and only if, $u
=\sum\limits_{g
\in G}
\a(g) e(g) \otimes x$ for some $x \in G$ and a map $\a : G \map \BC^*$ such that
$$
\g_x(g,h) =  \dfrac{\a(g)\a(h)}{\a(g h)}
$$
for any $g,h \in G$.
\end{prop}
\pf The result follows from direct computation. \qed\\

Let $\widehat{G}$ be the character group of $G$. Since $\g_1 = 1$, $\sum\limits_{g \in G}\a(g)
e(g) \otimes 1 \in \G^\w$ for any $\a \in
\widehat{G}$ by Proposition \ref{p0.1}. The map $\widehat{G} \map \G^\w$, $\a \mapsto
\sum\limits_{g
\in G}\a(g) e(g)
\otimes 1$ is an injective group homomorphism. In the sequel, we often identify $\widehat{G}$
with the image of $\widehat{G}$ under this map. On the other hand, the map
\begin{equation}\label{eq0.5}
D^\w(G)
\map
\BC G,\quad\sum\limits_{g,x \in G} \a(g,x) e(g) \otimes x \mapsto \sum\limits_{x\in
G}\a(1,x) x
\end{equation}
is an algebra map. Let $B^\w$ be the image of $\G^\w$ under this map.
One can easily see that
$$
B^\w = \{x \in G\,|\, \g_x \mbox{ is a 2-coboundary.}\}\,.
$$
\begin{lem} \label{l0.1}
Let $G$ be a finite group. Then $\widehat{G}$ is in the center of $D^\w(G)$.
Moreover, $\G^\w$ is a central extension
$$
\xymatrix{
1\ar[r]&\widehat{G}\ar[r]&\G^\w \ar[r] &B^\w\ar[r]& 1\,.
}
$$
For each $x \in
B^\w$,  let  $\g_x=\delta\tau_x$ for a 1-cochain $\tau_x:G \map \BC^*$.
The 2-cocycle $\b$ associated to this central extension is
given by
\begin{equation*}\label{eq0.7}
\b(x,y)(g) = \dfrac{\tau_x(g) \tau_y(x^{-1} g x)}{\tau_{xy}(g)} \theta_g(x,y)\,.
\end{equation*}
\end{lem}
\pf
For any $\a \in \widehat{G}$ and $e(h)\otimes x \in D^\w(G)$,
\begin{eqnarray*}
(e(h)\otimes x)\cdot(\sum_{g \in G} \a(g) e(g)\otimes 1) &=& \sum_{g\in
G}\a(g)e(h)e(x^{-1}g x)\otimes x \\
&=& \sum_{g\in
G}\a(xgx^{-1})e(h)e(g)\otimes x \\
&=& \sum_{g\in
G}\a(g)e(h)e(g)\otimes x \\
&=& (\sum_{g \in G} \a(g) e(g)\otimes 1) \cdot (e(h)\otimes x) \,.
\end{eqnarray*}
Therefore, the first statement follows.\\
 An element $u \in \ker \left(\G^\w\map
B^\w\right)$ if, and only if,
$u = \sum\limits_{g\in G}\a(g) e(g) \otimes 1$ and $\g_1=\delta \a$.  Since $\g_1 =1$,
$\a \in \widehat{G}$. Hence,
$$
\widehat{G} = \ker \left(\G^\w \map B^\w\right)\,.
$$
Since $\g_x= \delta\t_x$ for $x \in B^\w$,
$$
x \mapsto \sum_{g \in G} \t_x(g) e(g) \otimes x
$$
is a section of the map $\G^\w \map B^\w$. One can check directly that
$$
\left(\sum_{g \in G} \t_x(g) e(g) \otimes x\right)\cdot\left(\sum_{g \in G} \t_y(g) e(g)
\otimes y \right) = \sum_{g \in G} \t_x(g) \t_y(x^{-1} g x) \theta_g(x,y) e(g) \otimes xy.
$$
Hence, the formula for the 2-cocycle associated to the central extension follows. \qed
\begin{remark}
{\rm
In general, $B^\w$ and the extension $\G^\w$ depend on the individual cocycle.
However, both of them are independent of the representative of the cohomology class $[\w]$ when
$G$ is abelian.
}
\end{remark}

\begin{prop}\label{p1.2}
Let $G$ be a finite abelian group and $\w, \w'$  normalized 3-cocycles of $G$.
\begin{enumerate}
\item[\rm (i)] $\G^\w$ is abelian.
\item[\rm (ii)] If $\w$ and $\w'$ are cohomologous, then $B^\w=B^{\w'}$ and the
central extensions
$$
1 \map \widehat{G} \map \G^\w\map B^\w \map 1 \quad\mbox{and}\quad
1 \map \widehat{G} \map \G^{\w'} \map B^{\w'} \map 1
$$
are equivalent.
\end{enumerate}
\end{prop}
\pf (i) By Proposition \ref{p0.4}, $D^\w(G)_0$ is a self-dual Hopf algebra and
hence $\G^\w$ is isomorphic to the group of group-like elements of $D^\w(G)_0^*$.
Let $\s_1$, $\s_2$ be group-like elements of $D^\w(G)_0^*$. Then,
$\s_1$, $\s_2: D^\w(G) \map \BC$ are algebra maps. Let $V_1$, $V_2$ be the
1-dimensional representations associated to $\s_1$ and $\s_2$ respectively. Then,
$V_1 \otimes V_2$ and $V_2 \otimes V_1$ are the 1-dimensional representations associated
to $\s_1 \ast \s_2$ and $\s_2 \ast \s_1$ respectively. Since
$V_1 \otimes V_2 \cong V_2 \otimes V_1$ as $\H$ modules,
$\s_1 \ast \s_2 = \s_2 \ast \s_1$. Therefore, $\G^\w$ is abelian. \\
(ii) Let $\w' = \w \delta b$ for some normalized 2-cochain $b : G \times G \map \BC^*$.
Then
$$
\w'_g(x,y) = \w_g(x,y) \delta b_g (x,y)
$$
where $b_g(x) = b(x,g)/b(g,x)$. Hence, $\w'_g$ is a 2-coboundary if, and only if, $\w_g$ is a
2-coboundary. Therefore, $B^{\w'} = B^\w$. Since the map $\Theta : D^\w(G)_0 \map
D^{\w'}(G)_0$ defined in Proposition \ref{p0.4} is an Hopf algebra isomorphism,
$\Theta (\G^\w) = \G^{\w'}$. Moreover, $\Theta$ satisfies the commutative
diagram :
$$
\xymatrix{
1 \ar[r] & \widehat{G} \ar[r]\ar@{=}[d] & \G^\w\ar[d]_\Theta \ar[r]  &
B^\w\ar[r]\ar@{=}[d] & 1 \\
1 \ar[r] & \widehat{G}\ar[r] & \G^{\w'}\ar[r] & B^{\w'}\ar[r] & 1 \,.
}
$$
\qed  \\
\begin{cor}\label{c1.1}
Let $G$ be a finite group and $\w$ a normalized 3-cocycle.
Then the following statements are equivalent :
\begin{enumerate}
\item[\rm (i)] $D^\w(G)$ is spanned by $\G^\w$.
\item[\rm (ii)] $G$ is abelian and $B^\w=G$.
\item[\rm (iii)] $G$ is abelian and $\w  \in Z^3(G,\BC^*)_{ab}$
\item[\rm (iv)] $D^\w (G)$ is a commutative algebra.
\end{enumerate}
\end{cor}
\pf (i) $\Leftrightarrow$ (ii). It follows from Lemma \ref{l0.1} that
$$
|\G^\w|=|\widehat{G}|\cdot|B^\w|\,.
$$
Therefore, $D^\w(G)$ is spanned by $\G^\w$ if and only if $|\widehat{G}|=|G|$ and
$|B^\w|=|G|$ which is equivalent to the statement (ii).\\
The equivalence of (ii) and (iii) follows
directly from
the definition of $Z^3(G,\BC^*)_{ab}$ and
equation (\ref{multiplication}).
(ii) $\Rightarrow$ (iv). Since $G$ is abelian and
$B^\w=G$, $\w_g$ is a coboundary for any $g \in G$. In
particular, $\w_g(x,y) = \w_g(y,x)$ for any $x, y \in G$.
Therefore, the multiplication on $\H$ is commutative.\\
(iv) $\Rightarrow$ (ii)\, Assume $D^\w(G)$ is commutative. By the surjectivity of the
algebra map defined in (\ref{eq0.5}), $G$ is abelian and hence $\theta_g =
\g_g = \w_g$. Moreover, the commutative multiplication in $D^\w(G)$ implies that
$$
\w_g(x,y) = \w_g(y,x)
$$
for any $x,y,g \in G$. Therefore, $\w_g$ is a 2-coboundary for any $g \in G$.
 \qed\\

In the sequel, if $G$ is a finite abelian group we will denote by $\Lambda_G [\w]$
the cohomology class in $H^2(B^\w,\widehat{G})$ which corresponds to the central
extension
$$
\xymatrix{
1 \ar[r] & \widehat{G}\ar[r] & \G^{\w}\ar[r] & B^{\w}\ar[r] & 1 \,.
}
$$
In particular, if $[\w] \in H^3(G,\BC^*)_{ab}$, then $B^\w= G$ and so
$\Lambda_G[\w] \in H^2(G,\widehat{G})$. We  simply write $\Lambda$ instead of
$\Lambda_G$ when there is no ambiguity.  \\
\begin{remark}
{\rm If one of the conditions (i)--(iv) of Corollary \ref{c1.1} holds then we have
an isomorphism of Hopf algebras $\H_0 \cong \BC[\G^\w]$ }
\end{remark}

\begin{prop}\label{p2.7}
The map $\Lambda : H^3(G,\BC^*)_{ab} \map H^2(G,\widehat{G})$,
$[\w] \mapsto \Lambda[\w]$ is a group homomorphism.
\end{prop}
\pf Let $[\w] \in H^3(G,\BC^*)_{ab}$. Then for each $x \in G$, there is a normalized
1-cochain $\tau_x$ such that $\delta \tau_x = \w_x$. Then  $\Lambda[\w]$ can be
represented by $\b$ where $\b : G \times G \map \widehat{G}$ is given by
\begin{equation}\label{eq2.8}
\b(x,y)(z) = \frac{\tau_x(z)\tau_y(z)}{\tau_{x y}(z)} \w_z(x,y) \,.
\end{equation}
The result follows easily from this formula. \qed\\

Let $G$ be a finite abelian group and $\w \in Z^3(G,\BC^*)_{ab}$.
Denote by $T(\w)$
the set of all
normalized 2-cochains $\t$ on $G$ such that
$$
\w_x=\delta\t_x
$$
for $x \in G$ where $\t_x(y)=\t(x,y)$.
Pick any $\t \in T(\w)$. Denote by $\s(\a,x)$ the group-like element
$$
\sum_{g\in G} \a(g) \t_x(g)\, e(g)\otimes x
$$
for any $\a \in \widehat{G}$ and $x \in G$. We will write $\s_\t(\a,x)$ if we wish
to emphasize the dependence of $\s$ on $\t$. By Lemma
\ref{l0.1}, for any
$u
\in
\G^\w$, there exist unique $\a \in \widehat{G}$ and $x \in G$ such
that $u=\s(\a,x)$. It follows from Proposition \ref{p2.7} that
\begin{equation}\label{eq:sigma}
\s(\a,x)\s(\l,y)  =\s(\a\l \b(x,y),xy)
\end{equation}
for any $\a,\l \in \widehat{G}$ and $x,y \in G$ where $\b$ is given by equation
(\ref{eq2.8}). We  will simply write $\Lambda (\w)$ for
the 2-cocycle given in equation (\ref{eq2.8}) whenever the context is clear.
%%%%%
%
%
%\input tp.tex
\section{Tensor Products}\label{s3}
We consider here the structure of the tensor product of two quantum doubles.
Let $G$, $H$ be finite groups and $\w,\w'$ normalized 3-cocycles on $G$ and $H$ respectively
(with $\BC^*$ coefficients as always). Of course there are canonical projections
$G \times H \map G$, $G \times H \map H$, and we let  $\infl \w$, $\infl \w'$ denote the
3-cocycles in $Z^3(G\times H,\BC^*)$ obtained by inflating $\w$, $\w'$ using these projections.
Let $\zeta = (\infl \w)( \infl \w')$.
\begin{prop}\label{pTP.1}
There is a natural isomorphism of quasi-triangular quasi-Hopf algebras
$$
\xymatrix{
\iota : \H \otimes D^{\w'}(H) \ar[r]^-{\cong} & D^\zeta(G\times H)
}
$$
given by $\iota: e(g) \otimes x \otimes e(h) \otimes y \mapsto e(g,h) \otimes (x,y)$
for any $g,x \in G$ and $h,y \in H$.
\end{prop}
This means that $\iota$ is an isomorphism of bialgebras, and that it maps the Drinfeld
associators, $\cal R$-matrices and elements $\a$, $\b$ (as given in section \ref{s1}) to the
corresponding objects associated with $D^\zeta(G\times H)$. The calculations needed to prove
the proposition are essentially routine and we omit the proof.

The group-like elements in $D^\w(G) \otimes D^{\w'}(H)$ are of the form
$$
\left(\sum_{g \in G} \l_x(g) e(g) \otimes x \right) \otimes \left(\sum_{h \in H} \l'_y(h) e(h)
\otimes y \right)
$$
where $x \in B^\w$, $y \in B^{\w'}$ and $\l_x$ and $\l'_y$ are functions on $G$ and $H$
respectively such that $\delta\l_x = \g_x$, $\delta\l'_y =\g'_y$ and $\g$, $\g'$ are those
$\g$'s associated to $\w$ and $\w'$ respectively (cf.(\ref{eq0.002})). It is easy to see that
the group of all group-like elements in
$D^\w(G) \otimes D^{\w'}(H)$ is naturally isomorphic to $\G^\w(G) \times \G^{\w'}(H)$. As $\iota$
induces an isomorphism from the group of group-like elements of $D^\w(G) \otimes D^{\w'}(H)$ to
$\G^\zeta(G\times H)$,
$$
\G^\w(G) \times \G^{\w'}(H) \stackrel{\iota}{\cong} \G^\zeta(G \times H)\,.
$$

\begin{prop}\label{pTP.2}
Let $G$ be a finite group given as a direct product of groups $G=H \times K$ with
$|H|$ and $|K|$ coprime. Then for any normalized 3-cocycle $\w$ of $G$, there exist
normalized 3-cocycles $\eta$, $\eta'$ on $K$, $H$ respectively such that
$\H$ and $D^\eta(H) \otimes D^{\eta'}(K)$ are equivalent as quasi-triangular quasi-bialgebras.
\end{prop}
\pf It is well-known \cite{Brown} that $H^n(G,\BC^*) = \infl (H^n(H,\BC^*)) \infl (H^n(K,\BC^*))$
(written multiplicatively), and in particular, for any normalized 3-cocycle $\w$ of $G$, $\w$ is
cohomologous to  $\zeta=(\infl \eta) (\infl \eta')$ for some normalized 3-cocycles $\eta$,
$\eta'$ on $K$, $H$ respectively. Hence, $\H$ and $D^\eta(H) \otimes D^{\eta'}(K)$ are
equivalent as quasi-triangular quasi-bialgebras by Proposition \ref{pTP.1} and Remark
\ref{r0.1} (ii).\qed\\

We state formally what obtains when $G$ is abelian since we will need this case.
\begin{prop}
Let $G$ be a finite abelian group with decomposition $G=P_1 \times \cdots
\times P_k$ into the product of its Sylow
subgroups. Let
$\w \in Z^3(G,\BC^*)$. Then there are cocycles
$\eta_i \in Z^3(P_i, \BC^*)$ such that
$D^\w(G)$ and $D^{\eta_1}(P_1) \otimes \cdots \otimes D^{\eta_k}(P_k)$ are equivalent
as quasi-triangular quasi-bialgebras.
\end{prop}

Because of this result, many of the questions we address regarding $D^\w(G)$ ($G$ abelian) can
be reduced to the case that $G$ is a $p$-group for some prime $p$.
%%%%%%%
%
%
%\input section3.tex
\section{Eilenberg-Maclane Cohomology}\label{s4}
Let $G$ be a group. The bar complex of $G$, denoted
by $A_0(G)$, is the chain complex of abelian group
$$
\cdots \stackrel{\p}{\map} C_3 \stackrel{\p}{\map}  C_2
\stackrel{\p}{\map} C_1
\map \BZ
$$
where $C_n$ is the free abelian group generated by all
$n$-tuples $(x_1,\dots,x_n)$ of elements $x_i$ of $G$
and
$\p$ is a
$\BZ$-linear map defined by
$$
\p(x_1,\dots,x_n) = (x_2,\dots,x_n) +
\sum_{i=1}^{n-1}(-1)^i(x_1,\dots,x_{i-1},x_ix_{i+1},\dots,x_n)
+(-1)^n(x_1,\dots,x_{n-1})\,.
$$
We will call $(x_1,\dots,x_n) \in C_n$ a $n$-dimensional cell of
$A_0(G)$. For any abelian group $\Pi$,
$\Hom(A_0(G),\Pi)$ is a cochain complex. We will  denote by
$C^n(G,\Pi)$, $Z^n(G,\Pi)$, $B^n(G,\Pi)$ and
$H^n(G,\Pi)$
the dimension $n$-cochains, cocycles, coboundaries
and cohomology classes of
$\Hom(A_0(G),\Pi)$ respectively.

 For any two cells $(x_1,\dots,x_n)$ and $(y_1,\dots,y_m)$ in
$A_0(G)$, we can define a ``shuffle'' of the $n$  letters
$x_1,\dots,x_n$ through the
$m$ letters
$y_1,\dots,y_m$ to be the $n+m$-tuple in which the order of the $x$'s
and the order of $y$'s are preserved. The sign of the shuffle is the
sign of the permutation required to bring the shuffled letters back
to the standard shuffle $(x_1,\dots,x_n, y_1,\dots,y_m)$. Then, we
can define the ``star'' product of $(x_1,\dots,x_n)*
(y_1,\dots,y_m)$ to be the signed sum of the shuffles of the letters
$x$ through the letters $y$.

Let $G$ be an abelian group. In the paper \cite{Mac52}, the author
described a complex $A_1(G)$ in which the cells are symbols
$\sigma=[\a_1|\a_2|\cdots |\a_p]$, with each $\a_i$ a cell of
$A_0(G)$. The dimension of $\sigma$ is $p-1$ plus the sum of the
dimensions of the $\a_i$, and the boundary of $\sigma$ is
$$
\p \sigma = \sum_{i=1}^p(-1)^{\e_{i-1}}[\a_1|\cdots|\p
\a_i|\cdots|\a_p] + \sum_{i=1}^p
(-1)^{\e_i}[\a_1|\cdots|\a_i*\a_{i+1}|\cdots|\a_p]\,,
$$
where $\e_i = 1+ \dim[a_1|\cdots|\a_i]$. For any abelian $\Pi$,
$\Hom(A_1(G),\Pi)$ is a cochain complex. We  denote by
$C_{ab}^n(G,\Pi)$, $Z^n_{ab}(G,\Pi)$, $B^n_{ab}(G,\Pi)$ and
$H^n_{ab}(G,\Pi)$
the dimension $n$-cochains, cocycles, coboundaries
and cohomology classes of
$\Hom(A_1(G),\Pi)$ respectively. Until the end of this section, we will
tacitly assume that $\Pi$ is a divisible abelian. The underlying
product in $\Pi$ will be written additively.

Notice that $A_0(G)$ is a subcomplex of
$A_1(G)$. We have the exact sequence of complexes
$$
0 \map A_0(G) \map A_1(G) \map B(G) \map 0
$$
with $B(G)$ the quotient complex $A_1(G)/A_0(G)$.
Then
$$
0 \map \Hom(B(G),\Pi) \map
\Hom(A_1(G),\Pi)
\map
\Hom(A_0(G),\Pi)
\map 0 \,.
$$
is exact and we have the long exact sequence
$$
\cdots \stackrel{\d}{\map}
H^*(B(G),\Pi) \map H^*_{ab}(G,\Pi)
\map H^*(G,\Pi) \stackrel{\d}{\map}
H^{*+1}(B(G),\Pi)\map H^{*+1}_{ab}(G,\Pi)
\map\cdots
$$
where $H^n(B(G),\Pi)$ is the $n$-th cohomology group of the
cochain complex
$\Hom(B(G),\Pi)$.

The low dimensional cells in $B(G)$ are given in the following table :
\begin{center}
\begin{tabular}{c||c}
dimension & cells\\ \hline
1 & none \\
2 & none \\
3 & $[x|y]$ \\
4 & $[x,y|z]\,,\quad [x|y,z]$ \\
5 & $[x|y|z]\,,\quad [x,y,z|u]\,,\quad [u|x,y,z]\,,
\quad [x,y|u,v] $
\end{tabular}
\end{center}
Since $H^2_{ab}(G,\Pi) =\Ext_{\BZ}^1 (G,\Pi) =0$, we have the exact
sequence
\begin{equation}\label{eq:1}
0\map H^2(G,\Pi) \stackrel{\d}{\map}
H^3(B(G),\Pi) \map H^3_{ab}(G,\Pi) \map H^3(G,\Pi)
\stackrel{\d}{\map}  H^4(B(G),\Pi)\,.
\end{equation}
We are going to analyze the sequence (\ref{eq:1}) in some detail.
Let us  denote by $B_n(G)$ the group of $n$-chains of $B(G)$.
For any cocycle $b : B^3(G) \map \Pi$ and $x,y,z \in G$,
\begin{eqnarray}
b[\p(x,y)|z] &=& 0 \\
b[x|\p(y,z)] &=& 0 \,.
\end{eqnarray}
Hence, $b$ is a bicharacter of $G$. Conversely, it is also  easy
to see that any bicharacter on $G$ defines an element in
$H^3(B(G),\Pi)$ uniquely. Hence, we have the following lemma.
\begin{lem}
$H^3(B(G),\Pi) = \Hom(G\otimes G,\Pi)$.
\end{lem}
Let $c \in Z^2(G,\Pi)$ and $[c]$ the corresponding cohomology
class. There exists a natural isomorphism
$\theta : H^2(G,\Pi) \map \Hom(\Wedge^2 G,\Pi)$
defined by
$$
\theta ([c])(x,y) = c(x,y)-c(y,x)
$$
for any $x,y \in G$(cf.\cite{Brown}). On the other hand, the connecting
homomorphism \\ $\d : H^2(G,\Pi) \map H^3(B(G),\Pi)$ is given by
$$
\d([c])[x|y] = c((x)*(y))
$$
for $x,y \in G$. One can easily check that the following diagram
commutes :
\begin{equation}
\begin{CD}
H^2(G,\Pi) @>\d>> H^3(B(G),\Pi) \\
@V{\theta}V\wr V            @VV{\wr}V\\
\Hom(\Wedge^2 G, \Pi) @>>> \Hom(G \otimes G,\Pi)\,.
\end{CD}
\end{equation}
The bottom horizontal map in the diagram is the natural embedding.
The exact sequence (\ref{eq:1}) is thus equivalent to the following exact
sequence :
\begin{equation}\label{eq:2}
0\map \Hom(\Wedge^2 G, \Pi) \map  \Hom(G \otimes
G,\Pi) \map H^3_{ab}(G,\Pi) \map H^3(G,\Pi)
\stackrel{\d}{\map}  H^4(B(G),\Pi)\,.
\end{equation}

The 3-dimensional Eilenberg-Maclane cohomology
group $H^3_{ab}(G,\Pi)$ was explicitly described in \cite{Mac52}. The
elements in $Z^3_{ab}(G,\Pi)$ are pairs $(f,d)$ where $f \in
Z^3(G,\Pi)$ and $d \in C^2(G,\Pi)$ satisfying the following
conditions :
\begin{eqnarray}
 d(x y|z)-d(x|z)-d(y|z) +
f(x,y,z)-f(x,z,y)+f(z,x,y) &=& 0 \\
d(x|y z)-d(x|y)-d(x|z)
-f(x,y,z)+f(y,x,z)-f(y,z,x) &=& 0\,.
\end{eqnarray}

The 3-dimensional cocycle
$(f,d)$ is a coboundary if there exists $h \in
C^2(G,\Pi)$ such that
\begin{eqnarray}
f(x,y,z)&=&\d h (x,y,z)\\
d(x|y)&=&h(x,y)-h(y,x)\,.
\end{eqnarray}

Thus, the maps $\Hom(G\otimes G,\Pi) \map H_{ab}^3(G,\Pi)$ and
$H_{ab}^3(G,\Pi) \map H^3(G,\Pi)$ in (\ref{eq:2}) are given by
\begin{eqnarray}
  b \mapsto [(0,b)] &\mbox{and} & [(f,d)] \mapsto [f]
\end{eqnarray}
respectively, where $[x]$ means the cohomology class associated to
the cocycle $x$.\\
To any $(f,d) \in Z^3_{ab}(G,\Pi)$, one can assign the function
$t(x)=d(x|x)$, called its trace. Any trace is a
{\em quadratic function}--a function $t: G \map \Pi$ such that
\begin{enumerate}
\item[\rm (i)] $t(ax)=a^2t(x)$ for $a \in \BZ$, and
\item[\rm (ii)] $b_t(x,y)=t(x+y)-t(x)-t(y)$ defines a bilinear function on  $G$.
\end{enumerate}
\begin{thm}[Eilenberg-MacLane]\label{t3.2}
Let $Q(G,\Pi)$ consist of all quadratic
functions from $G$ to $\Pi$. The function assigning to each cocycle its trace
induces an isomorphism  
$$H^3_{ab}(G,\Pi)\stackrel{\cong}{\map}Q(G,\Pi)\,.$$
\end{thm}

A map $F : B^4(G)
\map \Pi$ is a cocycle if it satisfies the following:
\begin{eqnarray}
F[\p(x,y,z)|u] &=& 0\\
F[u|\p(x,y,z)] &=& 0\\
F[\p(x,y)|u,v] &=& F[x,y|\p(u,v)] \\
F[x|(y)*(z)] &=& - F[(x)*(y)|z]
\end{eqnarray}
$F$ is a coboundary if there is a function $d :
B^3(G) \map \Pi$ such that
\begin{eqnarray}
d[\p(x,y)|z] & = & F[x,y|z] \\
d[x|\p(y,z)] & = & F[x|y,z]\,.\\
\end{eqnarray}
Let $f \in Z^3(G,\Pi)$. The connecting homomorphism $\d :
H^3(G,\Pi) \map H^4(B(G),\Pi)$ in (\ref{eq:2}) is given by 
$\d[f]=[F]$ where
$F: B^4(G) \map \Pi$ is defined by
\begin{eqnarray}
F[x,y|z] &=&  -f((x,y)*(z)) \label{eq:3}\\
F[x|y,z] &=& f((x)*(y,z)) \label{eq:4}
\end{eqnarray}
Let $F\: B^4(G) \map \Pi$ be a cocycle. Define
\begin{equation}
\overline{F}(x,y,z) = F[x,y|z]-F[y,x|z]
\end{equation}
for $x,y,z \in G$.
\begin{lem}\label{l2}
The map $\Psi :[F] \mapsto \overline{F}$ defines a group homomorphism from
$H^4(B(G),\Pi)$ onto $\Hom(\Wedge^3 G,\Pi)$. The map splits if $G$ is
finitely generated.
\end{lem}
\pf Since $\p(x,y) =\p(y,x)$, one can check that $\Psi$ is a
well-defined. Let $F: B^4(G) \map \Pi$ be a
cocycle. Then, for any $x,y,u,v \in G$,
\begin{eqnarray*}
\overline{F}(\p(u,v),x,y) &=& F[\p(u,v),x|y]-F[x,\p(u,v)|y]\\
&=& -F[\p(u,v)|x,y]+F[\p(u,v)|y,x] \\
&=& -F[u,v|\p(x,y)]+F[u,v|\p(y,x)] \\
&=&0\,.
\end{eqnarray*}
Therefore, $\overline{F}$ is linear at the first entry. Clearly,
$$
\overline{F}(x,y,z)= -\overline{F}(y,x,z)
\quad\mbox{ and }\quad \overline{F}(x,x,z)=0 \,.
$$
In particular, $\overline{F}$ is also linear at the second entry.
Notice that
\begin{eqnarray*}
\overline{F}(x,y,z) &=& F[x,y|z]-F[y,x|z]\\
&=& F[x|z,y] - F[x|y,z]\\
&=& -(F[x,z|y]-F[z,x|y]) \\
&=& -\overline{F}(x,z,y)\,.
\end{eqnarray*}
Therefore, $\overline{F}$ is also linear at the third entry and hence
$\overline{F} \in \Hom(\Wedge^3 G,\Pi)$.

Note that the map $\psi^*: H^3(G,\Pi) \map \Hom(\Wedge^3G, \Pi)$ given by
\begin{equation}\label{phi*}
\psi^*[\w](x,y,z) = \sum_\s
\sgn(\s)\w(\s(x),\s(y),\s(z) )
\end{equation}
with $\s$ running through the permutations of $x,y,z$ is a surjective homomorphism (cf.
\cite{Brown}). It follows directly from (\ref{eq:3}) and (\ref{eq:4}) that the diagram
$$
\xymatrix{ H^3(G,\Pi)\ar[r]^-{\d} \ar[ddr]_{\psi^*} &  H^4(B(G),\Pi)\ar[dd]^{-\Psi} \\
 & \\
&\Hom(\Wedge^3G,\Pi)\\
}
$$
commutes.
Hence, $\Psi$ is surjective. Moreover, if
$G$ is finitely generated, the map $\psi^*$ has a linear section $s$. Then, $-\d s$ is a
linear section of $\Psi$. \qed\\

For any $\b \in Z^2_{ab}(G,\Hom(G,\Pi))$, define $F_\b : B^4(G) \map \Pi$ by
\begin{equation}\label{eq:5}
F_\b[x,y|z]= 0  \quad\mbox{and}\quad F_\b[x|y,z]=\b(y,z)(x)
\end{equation}
for any $x,y,z \in G$. One can easily check that $F_\b \in Z^4(B(G),\Pi)$. Moreover,
$\Xi : [\b] \map [F_\b]$ defines a linear map from $H^2_{ab}(G,\Hom(G,\Pi))$ to
$H^4(B(G),\Pi)$.
\begin{lem}
The sequence
$$
\xymatrix{
0 \ar[r]& H^2_{ab}(G,\Hom(G,\Pi)) \ar[r]^-{\Xi}& H^4(B(G),\Pi) \ar[r]^-{\Psi}& \Hom(\Wedge^3
G,\Pi) \ar[r] &0
}
$$
is exact. In particular, if $G$ is finitely generated, the sequence is split exact.
\end{lem}
\pf If $[F_\b]=0$, then there exists a function $d:B_3(G) \map \Pi$ such that
\begin{eqnarray*}
0&=& d[\p(x,y)|z]\\
\b(y,z)(x)&=& d[x|\p(y,z)]\,.
\end{eqnarray*}
Hence, $d \in C^1(G,\Hom(G,\Pi))$ and $\b=\delta d$. Therefore,
$\Xi$ is injective.\\
Obviously, $\Psi\Xi[\b] = 0$ for $\b \in
Z^2_{ab}(G,\Hom(G,\Pi))$. Let $[F] \in \ker \Psi$. Then, we have
$$
F[x,y|z] = F[y,x|z]
$$
for any $x,y,z \in G$. Since $\Pi$ is divisible, there exists a map $\t \in
C^2(G,\Pi)$ such that
$$
\t[\p(x,y)|z]=F[x,y|z]\,.
$$
Then the map $F':B^4(G) \map \Pi$ defined by
$$
F'[x,y|z] = \t[\p(x,y)|z] \quad\mbox{and}\quad F'[x|y,z]=\t[x|\p(y,z)]
$$
is a coboundary and hence $[F] = [F-F']$. Let $\b(y,z)(x) = F[x|y,z]-F'[x|y,z]$.
Then,
\begin{eqnarray*}
\b(y,z)(u)-\b(y,z)(u v)+\b(y,z)(v) & = & F[\p(u,v)|y,z] - \t(\p(u,v)|\p(y,z)) \\
&=& F[u,v|\p(y,z)]-\t(\p(u,v)|\p(y,z)) \\
&=& \t(\p(u,v)|\p(y,z))-\t(\p(u,v)|\p(y,z)) =0\,.
\end{eqnarray*}
Moreover, $\b(\p(x,y,z))(u)= F[u|\p(x,y,z)]-F'[u|\p(x,y,z)] =0$. Therefore, $\b \in
Z^2_{ab}(G,\Hom(G,\Pi))$ and $F_\b = F-F'$. Thus, $\Xi[\b]=[F]$. The second statement
follows directly from Lemma \ref{l2}. \qed
%%%%%%
%
%
%\input section4.tex
\section{The Kernel Of $\Lambda$}\label{s5}
Let us return to the exact sequence (\ref{eq:2}) with $\Pi=\BC^*$. We then have the exact
sequence
\begin{equation}\label{eq:4.1}
0\map \Hom(\Wedge^2 G, \BC^*) \map  \Hom(G \otimes
G,\BC^*) \map H^3_{ab}(G,\BC^*) \map H^3(G,\BC^*)
\stackrel{\d}{\map}  H^4(B(G),\BC^*)\,.
\end{equation}
Notice that for any $(f,d) \in Z^3_{ab}(G,\BC^*)$, $f \in Z^3(G,\BC^*)$ and
$f_z$ is a coboundary for $z \in G$ because
$$
f_z(x,y) = \frac{f(z,x,y)f(x,y,z)}{f(x,z,y)}=\frac{d(z|x)d(z|y)}{d(z|x y)}\,.
$$
Therefore, $[f] \in H^3(G,\BC^*)_{ab}$ and hence the image of the map
$H^3_{ab}(G,\BC^*) \map H^3(G,\BC^*)$ is contained in $H^3(G,\BC^*)_{ab}$.

Let $f \in Z^3(G,\BC^*)_{ab}$. Then $\delta[f]$ admits the representative
$F \in Z^4(B(G),\BC^*)$ defined by
$$
F[x,y|z]=f((x,y)\ast (z))^{-1} \quad\mbox{and}\quad F[x|y,z]=f((x)\ast (y,z))\,.
$$
On the other hand, $\Xi\Lambda[f]$ has the representative $F' \in Z^4(B(G),\BC^*)$
given by
$$
F'[x,y|z]=1\quad\mbox{and}\quad F'[x|y,z]=\frac{\t_y(x)\t_z(x)}{\t_{yz}(x)}f_x(y,z)\,.
$$
where $f_x(y,z) = \frac{\t_x(y)\t_x(z)}{\t_x(y z)}$. Let
$$
d[x|y]=\t_y(x)\,.
$$
Then $F'/F = \delta d$. Hence $F'$ and $F$ represent the same cohomology class in
$H^4(B(G),\BC^*)$ and the diagram
$$
\begin{CD}
H^3(G,\BC^*) @>\d>> H^4(B(G),\BC^*) \\
@A{incl}AA            @AA{\Xi}A\\
H^3(G,\BC^*)_{ab} @>\L>> H^2_{ab}(G,\widehat{G})
\end{CD}
$$
commutes. Moreover, we obtain an exact sequence
\begin{equation}\label{eq4.2}
\xymatrix{
 \Hom(G \otimes G,\BC^*) \ar[r]&
H^3_{ab}(G,\BC^*) \ar[r] & H^3(G,\BC^*)_{ab} \ar[r]^-{\Lambda} & H^2_{ab}(G,\widehat{G})\,.
}
\end{equation}
Therefore, $\ker \Lambda = \coker \left(\Hom(G \otimes G,\BC^*) \map
H^3_{ab}(G,\BC^*)\right)$. Since the map $\Hom(G \otimes G,\BC^*) \longrightarrow
H^3_{ab}(G,\BC^*)$ is given by $b \mapsto [(1,b)]$,
$$
\Im\left(\Hom(G \otimes G,\BC^*) \map
H^3_{ab}(G,\BC^*)\right)=
 \left\{[(1,b)] \in H^3_{ab}(G,\BC^*)\,|\, b \mbox{ is a bicharacter}\right\}\,.
$$
Let $K(G,\BC^*)$ be the group of quadratic forms which are the trace of some bicharacter.
Then, by Theorem \ref{t3.2}, we see that the following holds:
\begin{prop}
For any abelian group $G$, $\ker \Lambda  =
Q(G,\BC^*)/K(G,\BC^*)$.\qed
\end{prop}
\begin{lem}\label{l4.2}
Let $G$ be a finite abelian group.
\begin{enumerate}
  \item[\rm (i)] Let $G=H \oplus L$ for  subgroups $H,L$ of $G$. Then
  $t \in K(G,\BC^*)$ if, and only if, $t|_H \in K(H,\BC^*)$ and
  $t|_L \in K(L,\BC^*)$.
  \item[\rm (ii)] Let $H$ be the Sylow 2-subgroup of $G$. Then
  $$
  Q(G,\BC^*)/K(G,\BC^*) \cong Q(H,\BC^*)/K(H,\BC^*)
  $$
  as abelian groups.
  \item[\rm (iii)] $Q(G,\BC^*)^2 \C K(G,\BC^*)$.
\end{enumerate}
\end{lem}
\pf (i) The necessity of the statement is obvious. Let us assume
$t \in K(G,\BC^*)$ is such that $t|_H \in K(H,\BC^*)$ and
  $t|_L \in K(L,\BC^*)$. Let $b_H$ and $b_L$ be bicharacters of $H$ and
  $L$ respectively with
  $$
  t|_H(x) = b_H(x,x) \quad \mbox{and}\quad t|_L(y) = b_L(y,y)
  $$
  for $x \in H$ and $y \in L$. Since $t \in Q(G,\BC^*)$,
  $$c(u,v)=\frac{t(u v)}{t(u)t(v)}$$ defines a bicharacter on $G$.
  Consider the map $b:G \times G \map \BC^*$ given by
  $$
   b(x_1 y_1, x_2 y_2)=b_H(x_1,x_2)c(x_1,y_2) b_L(y_1,y_2)
  $$
  for any $x_1,x_2 \in H$ and $y_1,y_2 \in L$. One can check directly that
  $b$ is a bicharacter on $G$. Moreover,
  $$
  b(x y,x y)= b_H(x,x)c(x,y)b_L(y,y)=t(x y)\,.
  $$
  Therefore, $t \in K(G,\BC^*)$. \\
(ii) It suffices to show that $t \in K(G,\BC^*)$ if, and only if,
$t|_H \in K(H,\BC^*)$.
Let $L$ be the subgroup of odd order elements of $G$.
Then $G = H \oplus L$.
Since $|L|$ is odd, $K(L,\BC^*)=Q(L,\BC^*)$. Hence, by (i),
$t \in Q(G,\BC^*)$ if, and only if, $t|_H \in Q(H,\BC^*)$. \\
(iii) Since
$$
c(x,y) = \frac{t(x y)}{t(x)t(y)}
$$
is a bicharacter and $c(x,x)=t(x)^2$, the result follows.\qed

\begin{prop}\label{p4.3}
Let $G$ be a finite abelian group. Then $Q(G,\BC^*)/K(G,\BC^*) \cong
\Omega_2(G)$.In particular, $\ker \Lambda_G \cong \Omega_2(G)$.
\end{prop}
\pf By virtue of Lemma \ref{l4.2}(ii), we may assume that $G$ is a finite abelian 2-group.
Let us consider the case that $G$ is cyclic of order $2^n$ with the generator $y$.
Then
$$
1=t(1)=t(y^{2^n})=t(y)^{2^{2n}}\
$$
Since $y^{-1} = y^{2^n-1}$,
$$
t(y) = t(y^{2^n-1}) = t(y)^{2^{2n}-2^{n+1}+1} \,.
$$
Hence, $t(y)$ is a $2^{n+1}$th root of unity. Conversely, for any
$2^{n+1}$th root of unity $\xi$, $t(y^r)= \xi^{r^2}$ defines a quadratic
function. Hence $Q(G,\BC^*) \cong \BZ_{2^{n+1}}$. Obviously, $K(G,\BC^*)$
is a proper subgroup $Q(G,\BC^*)$ and $Q(G,\BC^*)^2 \C K(G,\BC^*)$
by Lemma \ref{l4.2}(iii). Since
$Q(G,\BC^*)^2$ is the largest proper subgroup of $Q(G,\BC^*)$,
$Q(G,\BC^*)^2 = K(Q,\BC^*)$. Hence,
$Q(G,\BC^*)/K(G,\BC^*) \cong \BZ_2 \cong \Omega_2(G)$. \\
Now, we consider the general case. Let $G = C_1 \oplus \cdots \oplus C_l$ where
$C_1, \dots , C_l$ are cyclic subgroups of $G$. Consider the map
$p:Q(G,\BC^*) \map Q(C_1,\BC^*)/K(C_1,\BC^*)\times \cdots \times
Q(C_l,\BC^*)/K(C_l,\BC^*) $ defined by
$$
p(t) = (t_1,\dots,t_l)
$$
where $t_i$ is the coset $t|_{C_i}K(C_i,\BC^*)$. Obviously, $p$ is an epimorphism.
By Lemma \ref{l4.2}(i), $\ker p = K(G,\BC^*)$. Hence,
$$
\frac{Q(G,\BC^*)}{K(G,\BC^*)} \cong \frac{Q(C_1,\BC^*)}{K(C_1,\BC^*)}\times \cdots \times
\frac{Q(C_l,\BC^*)}{K(C_l,\BC^*)}\,.
$$
As $\dfrac{Q(C_i,\BC^*)}{K(C_i,\BC^*)} \cong \BZ_2$,
$\dfrac{Q(G,\BC^*)}{K(G,\BC^*)} \cong \BZ_2^l \cong \Omega_2(G)$. \qed
\\

\begin{cor}\label{c4.4}
If $G$ is a finite abelian group of odd order, then $\Lambda$ is
injective.\qed
\end{cor}

\begin{cor}
Let $G$ be a finite abelian group and $H$ the Sylow 2-subgroup of
$G$. Let $H=\bigoplus_{i=1}^l C_i$ be a cyclic subgroup
decomposition of $H$. Let $p_i : H \map C_i$ be the
natural projection associated to the decomposition and
$p_i^*: H^3(C_i,\BC^*) \longrightarrow H^3(G,\BC^*)$ the associated inflation of $p_i$.
 Then
$$
\ker \Lambda_G = \sum_{i=1}^l p_i^*(\ker \Lambda_{C_i})
$$
\end{cor}
\pf Note that the exact sequence (\ref{eq4.2}) is natural in $G$.
Therefore, we have the commutative diagram
$$
\xymatrix{
 \Hom(G \otimes G,\BC^*) \ar[r]&
H^3_{ab}(G,\BC^*) \ar[r] & H^3(G,\BC^*)_{ab} \ar[r]^-{\Lambda_G} &
H^2_{ab}(G,\widehat{G}) \\
 \Hom(C_i \otimes C_i,\BC^*) \ar[r]\ar[u] &
H^3_{ab}(C_i,\BC^*) \ar[r] \ar[u] & H^3(C_i,\BC^*)_{ab} \ar[u]^-{p^*_i}
\ar[r]^-{\Lambda_{C_i}} &  H^2_{ab}(C_i,\widehat{C_i}) \ar[u]
}
$$
with exact rows. Therefore, $p_i^*(\ker \Lambda_{C_i}) \subseteq \ker \Lambda_G$
for $i=1,\dots,l$. Hence, $\sum_{i=1}^l p_i^*(\ker \Lambda_{C_i}) \subseteq \ker \Lambda_G$.
Since $p_i^*$ is injective and $\{\Im p_i^*\}$ is $\BZ$-linearly independent in
$H^3(G,\BC^*)_{ab}$,
$$
\sum_{i=1}^l p_i^*(\ker \Lambda_{C_i}) \cong \bigoplus_{i=1}^l\ker
\Lambda_{C_i} \cong \BZ_2^l\,.
$$
Therefore, the equality in the statement follows.\qed
%%%%%%
%
%
%\input section5.tex
\section{The Cohomology Group $H^3(G,\BC^*)_{ab}$}\label{s6}
\begin{defn}
{\rm
Let $G$ be an abelian group and $M$ be a left $G$-module. For $g \in G$, the group
homomorphism $m_g : C_{\bullet-1}(G) \map C_\bullet(G)$ given by
$y \mapsto y\ast g$ is a chain map.
We  denote by $D^\bullet_{g,M} \: H^\bullet(G,M) \map H^{\bullet-1}(G,M)$
the map induced by $m_g$ and define $D^0_{g,M}$ to be the trivial map.
We will also denote $H^\bullet(G,M)_{ab} = \bigcap_{g\in G}\ker D^\bullet_{g,M}$.
When $M$ is the trivial module $\BC^*$, $D_{g,M}$ will simply written as $D_g$. In this case,
one can easily see that $H^3(G,M)_{ab}$ coincides with $H^3(G,\BC^*)_{ab}$ defined previously.
}
\end{defn}

\begin{lem}\label{l5.2}
For any positive integer $n \ge 1$,
$$
H^n(G,\BC^*)_{ab} \cong H^{n+1}(G,\BZ)_{ab}
$$
and
$$
\frac{H^n(G,\BC^*)}{H^n(G,\BC^*)_{ab}} \cong \frac{H^{n+1}(G,\BZ)}{H^{n+1}(G,\BZ)_{ab}}\,.
$$
\end{lem}
\pf Recall that, for $n \ge 1$,  $H^n(G,\BC^*) \stackrel{\delta}{\cong}
H^{n+1}(G,\BZ)$ under the connecting maps of the long exact sequence
associated to the exact sequence of trivial $G$-modules
$$
\xymatrix{
0 \ar[r] & \BZ \ar[r] & \BC \ar[r]^-{exp} & \BC^* \ar[r] & 1 \,.
}
$$
Since the diagram
$$
\xymatrix{
H^n(G,\BC^*) \ar[r]^\delta \ar[d]_-{D^n_g} & H^{n+1}(G,\BZ)\ar[d]^-{D^{n+1}_{g,\BZ}} \\
H^{n-1}(G,\BC^*) \ar[r]^\delta & H^{n}(G,\BZ)
}
$$
commutes for $g \in G$, the result follows. \qed
\begin{remark}\label{r5.3}
{\rm
Let $G$ be an abelian group. For any $g \in G$,
the map $D_{g,\BZ} = \bigoplus_{n \ge 0} D^n_{g,\BZ}$
is a graded derivation on $\bigoplus\limits_{n \ge 0} H^n(G,\BZ)$ with respect
to the cup product. Hence, $\ker D_{g,\BZ}$ is a graded subring of $H^\bullet(G,\BZ)$.
Since $H^\bullet(G,\BZ)_{ab} =  \bigcap\limits_{g \in G} \ker D_{g,\BZ}$,
$H^\bullet(G,\BZ)_{ab}$ is a graded subring of $H^\bullet(G,\BZ)$. In particular,
if $G$ is a finite cyclic group, $H^\bullet(G,\BZ)$ is generated by $H^2(G,\BZ)$.
As $H^1(G,\BZ)=0$, $H^2(G,\BZ)= H^2(G, \BZ)_{ab}$. Therefore,
$H^\bullet(G,\BZ)_{ab}=H^\bullet(G,\BZ)$ in case $G$ is cyclic.
}
\end{remark}

\begin{lem}\label{l5.1}
 For any abelian group $G$, $H^2(G,\BC^*)_{ab}$ is trivial and
 $$
\frac{H^3(G,\BC^*)}{H^3(G,\BC^*)_{ab}} \cong \Hom(\Wedge^3 G,\BC^*)\,.
 $$
\end{lem}
\pf For $[f] \in H^2(G,\BC^*)_{ab}$, $f(x,y)/f(y,x)=1$. Hence, $f$ is a coboundary.
Recall that the map $\varphi^*\: H^3(G,\BC^*) \map \Hom(\Wedge^3 G, \BC^*)$ defined by
 equation (\ref{phi*}) is a split epimorphism. Moreover, for any
normalized 3-cocycle $\w$,
$$
\varphi*([\w]) (x,y,z) = \w_x(y,z)/\w_x(z,y)
$$
for any $x,y,z \in G$. Therefore, $\ker \varphi^* = H^3(G,\BC^*)_{ab}$. \qed\\

Let $G$ be a finite abelian group with a cyclic subgroup decomposition
$G = \bigoplus\limits_{i=1}^l C_i$\,. Let $\{p_i\}$ be the projections
of $G$ associated to the decomposition and $P^\bullet_i : H^\bullet(C_i,\BZ)
\map H^\bullet(G,\BZ)$  the induced graded ring homomorphism with respect to
cup product. Let $P^\bullet_{C_1,\dots,C_l;G} :
H^\bullet(C_1,\BZ) \otimes \cdots \otimes H^\bullet(C_l,\BZ)
\map H^\bullet(G,\BZ)$ be the graded ring monomorphism defined by
$$
P^\bullet_{C_1,\dots,C_l;G}(f_1 \otimes \cdots \otimes  f_l)
= P_1(f_1)\cup \cdots \cup P_l(f_l)\,.
$$

\begin{prop}\label{p5.4}
Let $G$ be a finite abelian group with a cyclic subgroup decomposition
$G = \bigoplus\limits_{i=1}^l C_i$\,.
Then
$$
H^4(G,\BZ)_{ab}= \Im P^4_{C_1,\dots,C_l;G}\,.
$$
In particular,
$$
H^4(G,\BZ)_{ab}=
\bigoplus_{n_1+ \cdots +n_l = 4} H^{n_1}(C_1, \BZ) \otimes \cdots \otimes
H^{n_l}(C_l, \BZ)\,.
$$
\end{prop}
\pf We proceed by induction on $l$. The case  $l=1$ follows from Remark \ref{r5.3}.
Assume $l >1$ and $H =\bigoplus\limits_{i=2}^{l} C_i$.
Then $\Im P^\bullet_{C_1,\dots,C_l;G} \C  \Im P^\bullet_{C_1,H;G}$. Moreover,
by the K$\ddot{\mbox{u}}$nneth theorem,
$$
\frac{H^4(G,\BZ)}{\Im P^4_{C_1,\dots,C_l;G}} \cong
\frac{\Im P^4_{C_1,H;G}}{\Im P^4_{C_1,\dots,C_l;G}}
\oplus \Tor(H^3(C_1,\BZ),H^2(H,\BZ))\,.
$$
Note that
$$
\frac{\Im P^4_{C_1,H;G}}{\Im P^4_{C_1,\dots,C_l;G}} \cong \frac{H^4(G,\BZ)}
{P^4_{C_2,\dots,C_l;H}}\,.
$$
By induction  and Lemma \ref{l5.1},
$$
\frac{\Im P^4_{C_1,H;G}}{\Im P^4_{C_1,\dots,C_l;G}} \cong \frac{H^4(G,\BZ)}
{H^4(G,\BZ)_{ab}} \cong \Wedge^3 G.
$$
Since
$$
\Tor(H^2(C_1,\BZ),H^3(H,\BZ)) \cong C_2 \otimes H^3(H,\BZ) \cong  C_2 \otimes
H^2(H,\BC^*) \cong C_2 \otimes  \Wedge^2 H \,,
$$
$$
\frac{H^4(G,\BZ)}{\Im P^4_{C_1,\dots,C_l;G}} \cong \left(C_2 \otimes  \Wedge^2 H \right)
\oplus \Wedge^3 H \cong \Wedge^3 G\,.
$$
Hence, by Lemma \ref{l5.1} and Lemma \ref{l5.2},
$$
\frac{H^4(G,\BZ)}{\Im P^4_{C_1,\dots,C_l;G}} \cong \frac{H^4(G,\BZ)}
{H^4(G,\BZ)_{ab}}\,.
$$
Thus, $|\Im P^4_{C_1,\dots,C_l;G}|=|H^4(G,\BZ)_{ab}|$.
Notice that $\Im P^\bullet_{C_1,\dots,C_l;G}$ is the graded subring of
$H^\bullet(G,\BZ)$ generated by $H^2(G,\BZ)$.
As $H^2(G,\BZ)= H^2(G,\BZ)_{ab}$, $\Im P^\bullet_{C_1,\dots,C_l;G}
\C H^\bullet(G,\BZ)_{ab}$. \qed
%%%%%
%
%
%\input section6.tex
\section{The Image of $\Lambda_G$}\label{s7}
Let $G$ be a finite abelian group and
$$
E: \quad 1 \map \widehat{G} \map \G \map G \map 1
$$
an abelian central extension of $G$ by $\widehat{G}$.  By applying the functor
$\widehat{?}=\Hom(?,\BC^*)$ to this exact sequence, we can obtain another
abelian central extension of $G$ by $\widehat{G}$, namely
$$
\e(E): \quad 1 \map \widehat{G} \map \widehat{\G} \map G \map 1
$$
where $G$ and $\widehat{\widehat{G}}$ are naturally identified. Obviously,
if $E_1$ and $E_2$ are equivalent abelian central extensions of $G$ by
$\widehat{G}$, and so are $\e(E_1)$ and $\e(E_2)$. Therefore, $\e$ induces
a map on $H^2_{ab}(G,\widehat{G})$. We will denote this map by the same symbol
$\e$.

Let $\b$ be a normalized abelian cocycle in $Z^2_{ab}(G,\widehat{G})$.
Then the set $\G=\widehat{G} \times G$
equipped with the multiplication
\begin{equation}\label{eq6.0-1}
(\a,x)\cdot (\l,y) = (\a\l \b(x,y),x y)
\end{equation}
is an abelian group. Moreover, a central extension of $G$ by $\widehat{G}$
associated to $[\b]$ is given by
\begin{equation}\label{eq6.0-2}
\xymatrix{
1 \ar[r] & \widehat{G} \ar[r]^i & \G \ar[r]^p & G \ar[r] & 1
}
\end{equation}
where $i(\a) = (\a,1)$ and $p(\a,x) = x$ for $\a \in \widehat{G}$ and $x \in G$.
Let
$\t_x$ be a normalized 1-cochain such that
\begin{equation}\label{eq6.0}
\delta \t_x (g,h)=\b(g,h)(x)
\end{equation}
for any $g,h \in G$.
For any $x \in G$, denoted by $\overline{x}$ be a fixed element in $\widehat{\G}$
such that $\widehat{i}(\overline{x})=x$ and $\overline{1}$ the identity element of
$\widehat{\G}$. Then, there exist a normalized cocycle $\b' \in Z^2(G,\widehat{G})$
such that
$$
\overline{x}\overline{y} = \widehat{p}(\b'(x,y))\overline{x}\cdot\overline{y}
$$
for any $x,y \in G$. Obviously, $\e[\b] = [\b']$.

Let $\t'_x(g)=\overline{x}(1,g)$ for any $g,x \in G$.
Then
$$
\b'(x,y)(g) = \frac{\t'_x(g)\t'_y(g)}{\t'_{x y}(g)}
$$
for $x,y,g \in G$. Since $\overline{x} \in \widehat{\G}$,
$$
\t'_x(g)\t'_x(h) = \t'_x(g h) \b(g,h)(x)\,.
$$
Therefore, $\t'_x = \t_x \l_x$ for some $\l_x \in \widehat{G}$.
 Define
\begin{equation}\label{eq6.1}
b'_1(x,y)(g) = \frac{\t_x(g)\t_y(g)}{\t_{x y}(g)} \quad \mbox{and}
\quad \l(x)(g) = \l_x(g)\,.
\end{equation}
Then $\b' = \b'_1 \delta \l$ and hence $\b'_1$ is a normalized 2-cocycle in
$Z^2_{ab}(G,\widehat{G})$ and
\begin{equation}\label{eq6.2}
\e([\b]) = [\b'_1] \,.
\end{equation}
By the formulae
(\ref{eq6.1}) and (\ref{eq6.2}), we obtain
\begin{prop}\label{p6.1}
The map $\e : H^2_{ab}(G,\widehat{G}) \map H^2_{ab}(G,\widehat{G})$ is a group
homomorphism. \qed
\end{prop}

Let us identify $\widehat{G}$ and $H^2(G,\BZ)$ and consider the map $\Upsilon :
\widehat{G}
\otimes
\widehat{G}
\map H^2_{ab}(G,\widehat{G})$, which assigns to  $\a\otimes \l$ the
cohomology class
$[\b]$ with
$$
\b(x,y)(g) = \a(g)^{\tilde{\l} (x,y)}
$$
where
$\tilde{\l}$ is a normalized 2-cocycle in $Z^2(G,\BZ)$ corresponding to
$\l$. Since $G$ is a finite, one can easily see that $\Upsilon$ is an
isomorphism.
\begin{lem}\label{l6.2}
Let $T :\widehat{G} \otimes \widehat{G} \map \widehat{G} \otimes \widehat{G}$ be the
transposition automorphism. Then,
$$
\Upsilon T  = \e \Upsilon\,.
$$
\end{lem}
\pf Let $\a_1, \a_2 \in \widehat{G}$ and $\tilde{\a_i}= \delta f_i$
for some $f_i : G \map \BC$ such that $\Im \delta f_i \C \BZ$ and $f_i(1)=0$.
Then $\Upsilon(\a_1 \otimes \a_2)$ is the cohomology class represented by $\b$ given
by
$$
\b(x,y)(g) = \a_1(g)^{f_2(x)+f_2(y) - f_2(x y)}\,.
$$
Let
$$
\t_g(x) = \a_1(g)^{f_2(x)}\a_2(x)^{f_1(g)} \,.
$$
Then $\b(x,y)(g) = \delta \t_g (x,y)$.
By formulae (\ref{eq6.1}) and (\ref{eq6.2}),
$\e [\b]=[\b']$ where $\b'$ is given by
$$
\b'(x,y)(g) = \frac{\t_x(g) \t_y(g)}{\t_{x y}(g)}
= \a_2(g)^{f_1(x)+f_1(y) - f_1(x y)}\,.
$$
Hence, $[\b'] = \Upsilon (\a_2 \otimes \a_1)$. \qed
\begin{lem}\label{l6.3}
For any normalized $\b \in Z^2_{ab}(G,\widehat{G})$,
$\xi(\b)(x,y,z) = \b(x,y)(z)$ defines a 3-cocycle in $Z^3(G,\BC^*)$ and
$[\xi(\b)] \in H^3(G,\BC^*)_{ab}$. Moreover,
$\Lambda_G[\xi(\b)] = [\b]\e[\b]$.
\end{lem}
\pf For any normalized $\b \in Z^2_{ab}(G,\widehat{G})$,
one can easily check that $\w = \xi(\b)$ satisfies the 3-cocycle identity.
For any $z \in G$, let $\t_z : G \map \BC^*$ such that
$\b(x,y)(z) = \delta \t_z (x,y)$. Then
$$
\w_z(x,y) = \delta \t_z (x,y)
$$
and so $[\w] \in H^2(G,\BC^*)_{ab}$. Moreover, $\Lambda_G[\w]$ can be
represented by the 2-cocycle $\b_1$ given by
$$
\b_1(x,y)(z) = \frac{\t_z(x) \t_z(y)}{\t_z(xy)}
\frac{\t_x(z) \t_y(z)}{\t_{xy}(z)}\,.
$$
However, from the foregoing we have $[\b_1]=[\b]\e[\b]$. \qed

\begin{thm}\label{t6.4}
Let $G$ be a finite abelian group. Then
$$
\Im \Lambda_G = \Im S
$$
where $S : H^2(G,\widehat{G}) \map  H^2(G,\widehat{G})$ is defined by
$S[\b] = [\b]\e[\b]$ for $[\b] \in  H^2_{ab}(G,\widehat{G})$.
\end{thm}
\pf It follows from Lemma \ref{l6.3} that $\Im \Lambda_G \C \Im S$.
By virtue of Lemma \ref{l6.2}, $\Im S \cong \Im (id + T)$.
Let $G=C_1 \oplus \cdots \oplus C_l$ be a cyclic subgroup decomposition
of $G$. As $\widehat{G} \cong G$, one can  easily see that
$$
\Im (id + T) \cong \left(\bigoplus_{i < j} C_i \otimes C_j\right) \oplus
\left(\bigoplus_{i=1}^l \frac{C_i}{(C_i)_2}\right)\,.
$$
By Proposition \ref{p5.4} and Proposition \ref{p4.3}
$$
H^3(G,\BC^*)_{ab} \cong \left(\bigoplus_{i < j} C_i \otimes C_j\right) \oplus
\left(\bigoplus_{i=1}^l C_i\right)
\quad \mbox{and}\quad \ker \Lambda_G \cong \Omega_2(G) \,.
$$
Therefore, $|\Im \Lambda_G|=|\Im S|$ and hence $\Im \Lambda_G = \Im S$.\qed
\begin{thm}\label{t6.5}
Let $G$ be a finite abelian 2-group and $\w$ a normalized 3-cocycle such that
$[\w] \in H^2(G,\BC^*)_{ab}$. Then $\G^\w(G)$ is a direct sum of
an {\rm even} number of
cyclic subgroups.
\end{thm}
\pf By Theorem \ref{t6.4}, $\Lambda_G[\w]=[\b] \e [\b]$ for some
normalized $\b \in Z^2_{ab}(G,\widehat{G})$. For $z \in G$, let
$\t_z :G \map \BC^*$ such that
$$
\b(x,y)(z) = \delta \t_z (x,y) \,.
$$
Then, $\Lambda_G[\w]$ contains the normalized 2-cocycle $\s$ given by
$$
\s(x,y)(z) = \frac{\t_z(x)\t_z(y)}{\t_z(x y)} \frac{\t_x(z)\t_y(z)}
{\t_{x y}(z)}\,.
$$
Consider the central extension associated to $\s$
\begin{equation}
\xymatrix{
1 \ar[r] & \widehat{G} \ar[r]^i & \G \ar[r]^p & G \ar[r] & 1
}
\end{equation}
as defined in (\ref{eq6.0-2}) where $\G$ the group with the underlying set
$\widehat{G} \times G$ endowed with the multiplication
$$
(\a,x)\cdot (\l,y) = (\a\l \s(x,y),x y)\,.
$$
To show that $\G$ is a direct sum of an even number of cyclic subgroups, it is
enough to show that $\G'=p^{-1}(\Omega_2(G))/i(\widehat{G}^2)$
is a direct sum of an even number of
cyclic subgroups. For this, it suffices to show  that $\G'$
admits a nonsingular alternating bilinear form valued in $\BC^*$ \cite{Wa63}.
Define $b : \G' \otimes \G' \map \BC^*$ by
$$
b((\a\widehat{G}^2,x),(\l\widehat{G}^2,y)) = \frac{\a(y)}{\l(x)} \frac{\t_y(x)}{\t_x(y)}\,.
$$
for $\a,\l \in \widehat{G}$ and $x,y \in \Omega_2(G)$. One can easily see that $b$ is
well-defined. The linearity of $b$ follows from
the fact that
$$
\b(x,y)=\frac{\t_x \t_y}{\t_{xy}}
$$
is a character for any $x,y \in G$. The routine verification that $b$ is also
non-degenerate and alternating will be left to the reader. \qed\\
\begin{remark}\label{r6.6}
{\rm
We summarize some of our results.
\item[(i)] We have an exact sequence
$$
1 \map \Hom(\Wedge^2G,\BC^*) \map \Hom(G \otimes G,\BC^*) \map
H^3_{ab}(G,\BC^*) \map H^3(G,\BC^*)_{ab} \stackrel{\Lambda_G}{\map} \Im S
\map 1
$$
where $S$ is as in Theorem \ref{t6.4}.
\item[(ii)] If $|G|$ is odd then $\Lambda$ induces an isomorphism
$$
\xymatrix{
H^3(G,\BC^*)_{ab} \ar[r]^\Lambda_\cong & H^2_{ab}(G,\widehat{G})^{+}
}
$$
where $H^2_{ab}(G,\widehat{G})^{+}$ is the group of $\e$-invariants. (In this case, $\Lambda$
is an injection (Cor. \ref{c4.4}) and $\Im S= H^2_{ab}(G,\widehat{G})^{+}$ as follows from
Theorem \ref{t6.4} and the fact that $\e$ is an involution on the group
$H^2_{ab}(G,\widehat{G})$.) }
\end{remark}
%%%%%%%%%%
%
%
%\input section7.tex
\section{The Monoidal Category $\M{\H}$}\label{s8}
Let $G$ be a finite abelian group and $\w \in Z^3(G,\BC^*)_{ab}$. Adopt the notation
introduced in section \ref{s2}, and pick $\t \in T(\w)$. Then,
$$
\G^\w =\{\s(\a,x)\,|\,\a \in \widehat{G}, x \in G \}
$$
where $\s(\a,x)$ is given by (\ref{eq:sigma}).
By Proposition \ref{p0.4}, for any $x \in G$ and $\a \in \widehat{G}$, the element
$$
\chi(\a,x) = \sum_{g \in G} \a(g)\t_x(g)f_{g,x}
$$
is an algebra map from $\H$ to $\BC$ where $\{f_{x,g}\}$ is the dual basis of
$\{e(g)\otimes x\}$. Denote by
$V(\a,x)$ the irreducible representation of $\H$ associated to the
algebra  map $\chi(\a,x)$. Then,
$$
S = \left\{V(\a,x)\left|\a \in \widehat{G}, x \in G\right\}\right.
$$
is a complete set of irreducible representations of $\H$ up to isomorphism and
$e(h)\otimes y$ acts on $V(\a,x)$ as  scalar multiplication by
$\a(y)\t_x(y) \delta_{h,x}$. Hence, as reflected in equation (\ref{eq0.01}), the
associativity map
\begin{equation}\label{eq7.01}
a_{V(\a,x),V(\l,y),V(\mu,z)} :
(V(\a,x) \otimes V(\l,y)) \otimes V(\mu,z)  \map
V(\a,x) \otimes (V(\l,y) \otimes V(\mu,z))
\end{equation}
in $\M{\H}$ is the scalar
$\w(x,y,z)^{-1}$. Notice that the group structure
on
$S$ induced by the tensor product of $\M{\H}$ is
isomorphic to $\G^\w$ with the isomorphism
$\s(\a,x)
\mapsto V(\a,x)$. Hence, by the results of \cite{TaYa98}, we have the
following  lemma.
\begin{lem}\label{l7.1}
Let $G$ be a finite abelian group and $\w$ be a normalized 3-cocycle such
that $[\w] \in H^3(G,\BC^*)_{ab}$. Then, $\M{\H}$ is monoidally equivalent
to $\M{\H_0}$  if, and only if,
$\infl [\w]$ is trivial,   where $\infl: H^3(G,\BC^*) \map
H^3(\G^\w,\BC^*)$ is  inflation along the projection $\G^\w \map G$.
\qed
\end{lem}
In fact, we have
\begin{lem}\label{r7.1}
Suppose that $G$ is finite and abelian and that $\infl [\w]$ is trivial. Then
$\H$ and $\H_0$ are gauge equivalent as quasi-bialgebras.
\end{lem}
\pf Since  $\infl [\w]$ is trivial,
there is a
 normalized 3-cochain $f$ on $\G^\w$ such that
\begin{equation}\label{eq7.1}
\infl(\w) = \delta f.
\end{equation}
Let $\t \in T(\w)$ and let
$E_{\s(\a,x)}=\frac{1}{|G|}e(x)\otimes \sum_{g\in G}\frac{g}{\a(g)\t_x(g)}$. One can easily
see that $E_u E_v = \delta_{u,v}E_v$ and $\Delta E_u = \sum_{s,t \in \G^\w} E_s \otimes E_t$.
Set
$$
F = \sum_{u,v \in \G^\w}f(u,v)^{-1} E_u \otimes E_v \,.
$$
One can easily check that $F$ is invertible with
$$
F^{-1} = \sum_{u,v \in \G^\w}f(u,v)E_u \otimes E_v \,.
$$
As $f$ is normalized,
$$
(\varepsilon \otimes id) F = 1_\H = (id \otimes \varepsilon)(F) \,.
$$
Moreover, by equation (\ref{eq7.1}),
$$
\Phi = F_{23} (id \otimes \Delta)(F)
\Phi_0 (\Delta \otimes id)(F^{-1})F_{12}^{-1}\,.
$$
where $\Phi_0 = 1_\H \otimes 1_\H \otimes 1_\H$.
Since $\H_0$ is a commutative, $\H$ is a twist of $\H_0$ by
$F$, that is, $\H$ and $\H_0$ are gauge equivalent.\qed

\begin{lem}\label{l7.2}
$\infl [\w] \in \ker
\Lambda_{\G^\w}$\,.
\end{lem}
\pf Set $\overline{\w}=\infl \w$.
Define $d : \G^\w \times \G^\w \map \BC^*$ by
$$
d[\s(\a,x)|\s(\l,y)]= \frac{1}{\a(y)\t_x(y)}\,.
$$
One can check directly that $(\overline{\w},d) \in Z^3_{ab}(\G^\w,\BC^*)$.
By the exact sequence (\ref{eq4.2}), we have the commutative diagram
$$
\xymatrix{
H^3_{ab}(G,\BC^*) \ar[r]\ar[d] & H^3(G,\BC^*)_{ab} \ar[r]^-{\Lambda}\ar[d]^\infl &
H^2_{ab}(G,\widehat{G})\ar[d]\\
H^3_{ab}(\G^\w,\BC^*) \ar[r] & H^3(\G^\w,\BC^*)_{ab}
\ar[r]^-{\Lambda}&
H^2_{ab}(\G^\w,\widehat{\G^\w})
}
$$
with the rows exact, where the vertical homomorphisms are induced by the
epimorphism $\G^\w \map G$. As
$$
[(\overline{\w},d)] \mapsto [\overline{\w}] = \infl [\w]\,,
$$
it follows from the diagram that $\infl[\w] \in \ker \Lambda_{\G^\w}$. \qed

\begin{thm}\label{t7.3}
Let $G$ be a finite group of odd order and $\w_1, \w_2$ normalized 3-cocycles such
that $[\w_1],[\w_2] \in H^3(G,\BC^*)_{ab}$. Then,
$\M{\Hone}$ and $\M{\Htwo}$
are monoidally equivalent if, and only if, $\G^{\w_1} \cong \G^{\w_2}$ as groups.
\end{thm}
\pf Since the group structure, induced by the underlying tensor product, on
the isomorphism classes of irreducible representations of $\Hi$ is
isomorphic to $\G^{\w_i}$, $\G^{\w_1} \cong \G^{\w_2}$ as
groups if $\M{\Hone}$ and $\M{\Htwo}$ are monoidal equivalent. Conversely,
assume that $\G^{\w_1} \cong \G^{\w_2}$ as groups. Then,
$\BC\G^{\w_1}$ and $\BC\G^{\w_2}$ are isomorphic as Hopf algebras.
By Lemma \ref{l7.2}, $\infl[\w_i] \in \ker \Lambda_{\G^{\w_i}}$.
As $G$ is of odd order, $|\G^{\w_i}|$ is odd, and it follows from Corollary
\ref{c4.4} that $\ker \Lambda_{\G^{\w_i}}$ is trivial. Hence,
by Lemma \ref{r7.1}, $\Hi$ is equivalent to $\BC\G^{\w_i}$ and
hence $\M{\Hi}$ is monoidally equivalent to the tensor category
$\M{\BC\G^{\w_i}}$ with the usual associativity constraint. \qed

\begin{example}\label{ex7.4}
{\rm
Theorem \ref{t7.3} is not true for groups of even order. For instance, take
$G = \BZ_2$. By Proposition \ref{p4.3}, $\ker \Lambda_G \cong \BZ_2$ and
hence $\ker \Lambda_G = H^3(G,\BC^*)$. Let
$\w_1$  be a normalized 3-cocycle on $G$ whose cohomology class
is nontrivial and let $\w_0=1$.   As $\Lambda_G [\w_i]$ is trivial
$(i=0,1)$,
$\G^{\w_i} \cong \BZ_2 \times \BZ_2$ as groups. Therefore,  inflation
of  cohomology along  $\G^{\w_i} \map G$ is injective for $i=0,1$. Thus,
$\infl [\w_1]$ is non-trivial, while  $\infl[\w_0]$ is obviously trivial. Thanks to
lemma \ref{l7.1}, if $\M{\BC[\BZ_2\times \BZ_2 ]}$  equipped with the usual
associativity constraint then $\M{D^{\w_0}(G)}$ is monoidally equivalent to
$\M{\BC[\BZ_2\times\BZ_2]}$ but $\M{D^{\w_1}(G)}$ is not. Note that $\BC[\BZ_4]$
and $\BC[\BZ_2 \times \BZ_2]$ are the only 4-dimensional semisimple Hopf
algebras up to isomorphism. As the fusion rule of $\M{D^{\w_1}(G)}$ is
isomorphic to $\BZ_2 \times \BZ_2$,  $\M{D^{\w_1}(G)}$ and $\M{\BC[\BZ_4]}$
are not monoidally equivalent. Hence, $D^{\w_1}(G)$ cannot be obtained from
any Hopf algebra by a twist.
}
\end{example}
\begin{remark}
{\rm
In the paper \cite{DPR90}, p92, the authors asked whether $D^\w(G)$ can be
obtained by twisting a Hopf algebra. Theorem \ref{t7.3} gives an
affirmative answer to the question under the  conditions stated in
the theorem. However, in general,  Example \ref{ex7.4} gives a negative
answer to the question.
}
\end{remark}
%%%%%%%%%%
%\input section8.tex
\section{Gauge Equivalence and Quadratic Forms}\label{s9}
Let $G$ be a finite abelian group, $\w \in
Z^3(G,\BC^*)_{ab}$ and $\t \in T(\w)$.
Then, as shown in the
proof of Lemma \ref{l7.2}, $(\infl \w^{-1},d_\w) \in
Z^3_{ab}(\G^\w(G),\BC^*)$ where $d_\w: \G^\w(G) \times \G^\w(G) \map \BC^*$
is defined by
$$
d_\w[\s(\a,x)|\s(\l,y)] = \a(y)\t_x(y)\,.
$$
Define $q_\w : \G^\w(G) \map \BC^*$ to be the trace of the abelian 3-cocycle $(\infl \w^{-1},d_\w)$. Obviously,
$q_\w$ is a quadratic map given by
\begin{equation}\label{eq8.1}
q_\w(\s(\a,x)) = \a(x)\t_x(x)\,.
\end{equation}

Denote by $b_\w$ the symmetric bicharacter on $\G^\w(G)$ associated to
$q_\w$, namely $b_\w(x,y)=\frac{q_\w(x y)}{q_\w(x)q_\w(y)}$.
One can easily show that $b_\w$ is given by
$$
b_\w(\s(\a,x),\s(\l,y))=\a(y)\l(x)\t_x(y)\t_y(x)
$$
for any $\s(\a,x),\s(\l,y) \in \G^\w(G)$.
Obviously, $b_\w$ is non-degenerate on $\G^\w(G)$. Let $\<\bullet,\bullet\>_\w$
be the $\BC$-linear extension of $b_\w$ on $\H$.
\begin{lem}\label{l8.1}
With the previous notation, the map $\varphi_\w\:\H_0 \map \H_0^*$,
$\varphi_\w\:x \mapsto \<x,\bullet\>_\w$ is identical to the Hopf algebra
isomorphism $\varphi: \H_0 \map \H_0^*$ defined in (\ref{eq0.4}). In particular,
$\<\bullet,\bullet\>_\w$ is a non-degenerate symmetric bilinear form on $\H$. Moreover,
for any $u, v\in \G^\w$,
\begin{equation}\label{eq8.04}
(\varphi_\w(u) \otimes \varphi_\w(v))R=d_\w[v|u]
\end{equation}
and
\begin{equation}\label{eq8.05}
q_\w(u) = (\varphi_\w(u) \otimes \varphi_\w(u))R
\end{equation}
where $R$ is the $\cal R$-matrix of $\H$.
\end{lem}
\pf To show that $\varphi_\w = \varphi$, it suffices to show that
$\varphi_\w(u) = \varphi(u)$ for $u\in \G^\w$. Hence, it is enough to show that
$$
\varphi(u)(v) = \varphi_\w(u)(v) = b_\w(u,v)
$$
for any $u,v \in \G^\w$. However, it follows from  (\ref{eq0.4})
 that
$$
\varphi(\s(\a,x))(\s(\l,y)) = \a(y)\l(x)\t_x(y)\t_y(x) = b_\w(\s(\a,x),\s(\l,y))\,.
$$
Now, $\varphi_\w(\s(\a,x))(e(g)\otimes y) = \a(y)\t_x(y)\delta_{g,x}$.
Thus,
$$
(\varphi_\w(\s(\a,x)) \otimes \varphi_\w(\s(\l,y)))R = \l(x)t_y(x) =
d_\w[\s(\l,y)|\s(\a,x)]
$$
and hence
$$
(\varphi_\w(u) \otimes \varphi_\w(u)R=d_\w[u|u]= q_\w(u)
$$
for $u \in \G^\w$. \qed

\begin{lem}\label{l8.2}
Let $\w_1$ and $\w_2$ be normalized 3-cocycles in $Z^3(G_1,\BC^*)_{ab}$ and
$Z^3(G_2,\BC^*)_{ab}$ respectively. If
$j:\DG{1}_0 \map \DG{2}_0$ is a bialgebra
isomorphism, then
$j^t: \DG{2}_0 \map \DG{1}_0$, defined by
$$
\< j^t(u),v \>_{\w_1} = \< u , j(v) \>_{\w_2}
$$
for any $u \in \DG{2}$ and $v \in \DG{1}$, is also a bialgebra isomorphism.
\end{lem}
\pf Consider the map $j':\DG{2}_0 \map \DG{1}(G)_0$ given by
$$
j' = \varphi^{-1}_{\w_1} j^* \varphi_{\w_2}\,.
$$
Since $j \: \DG{1}_0 \map \DG{2}_0$ is a bialgebra isomorphism, and so is
$j^* \: \DG{2}_0^* \map \DG{1}_0^*$. It follows from Lemma \ref{l8.1} that $j'$ is
also a bialgebra isomorphism. Moreover, for any $u \in \DG{2}$ and $v \in \DG{1}$,
$$
\< j'(u),v \>_{\w_1} = \left(\varphi_{\w_1}j'(u)\right)(v)=\left(j^*\varphi_{\w_2}(u)\right)(v)=
\varphi_{\w_2}(j(u))(v) =
 \< u , j(v) \>_{\w_2}\,.
$$
Hence, $j'=j^t$. \qed
\begin{defn}
{\rm Let $\G_1$, $\G_2$ be abelian groups, and $b_1$, $b_2$
bicharacters on $\G_1$ and
$\G_2$ respectively. The pairs $(\G_1,b_1)$ and $(\G_2,b_2)$ are said to be
{\em equivalent}
if there is a group isomorphism $j:\G_1 \map \G_2$ such that
$$
b_2(j(x),j(y)) = b_1(x,y)
$$
for $x,y \in \G_1$. Similarly, let $q_1 \:\G_1 \map \BC^*$ and $q_2\:\G_2 \map \BC^*$
be quadratic maps. The pairs $(\G_1,q_1)$ and $(\G_2,q_2)$, which we call {\em quadratic spaces},
are said to be {\em equivalent}
if there is a group isomorphism $j:\G_1 \map \G_2$ such that
$$
q_2 j = q_1\,.
$$
The {\em orthogonal sum} of two quadratic spaces $(\G_1,q_1)$ and $(\G_2,q_2)$ is the quadratic
space $(\G_1 \times \G_2, q)$ where $q(x,y)=q_1(x)q_2(y)$.
}
\end{defn}
\begin{thm}\label{t8.4}
Let $G_1$, $G_2$ be a finite abelian groups and $\w_1$, $\w_2$  normalized 3-cocycles in
$Z^3(G_1,\BC^*)_{ab}$ and $Z^3(G_2,\BC^*)_{ab}$ respectively.
Then, $\DG{1}$ and $\DG{2}$ are gauge equivalent  if, and only if,
$(\G^{\w_1},q_{\w_1})$ and $(\G^{\w_2},q_{\w_2})$ are equivalent.
\end{thm}
\pf Suppose that $\DG{1}$ and $\DG{2}$ are equivalent as quasi-triangular quasi
Hopf algebras. Then there exists a gauge transform $F \in \DG{2}\otimes \DG{2}$ such
that $\DG{1}$ and $\DG{2}_F$ are isomorphic as quasi-triangular quasi-bialgebras.
Let $j:\DG{1}\map \DG{2}_F$ be such an isomorphism. Then,
\begin{equation}\label{eq8.2}
(j\otimes j)R_1 = F_{21}R_2F^{-1}
\end{equation}
where $R_1$ and $R_2$ are the $\cal R$-matrices of $\DG{1}$ and $\DG{2}$ respectively.
Since $\Delta_2=(\Delta_2)_F$,
$\DG{2}_F$ and $\DG{2}_0$ are identical as bialgebras. By Lemma \ref{l8.2},
$j^t : \DG{2}_0 \longrightarrow \DG{1}_0$  is a  bialgebra isomorphism. In particular,
$j^t(\G^{\w_2})=\G^{\w_1}$.
Let $\{E_u\}_{u\in \G^{\w_2}}$ be the dual basis of $\G^{\w_2}$ with respect to the
pairing $\<\bullet,\bullet\>_{\w_2}$. Let $F = \sum_{u,v \in \G^{\w_2}}f(u,v) E_u \otimes E_v$.
For any  $u \in \G^{\w_2}$, apply $\varphi_{\w_2}(u) \otimes \varphi_{\w_2}(u)$
to equation (\ref{eq8.2}). The left hand side of the equation becomes
$$
(\varphi_{\w_2}(u) \otimes \varphi_{\w_2}(u))(j \otimes j)R_1=
(\varphi_{\w_1}(j^tu) \otimes \varphi_{\w_1}(j^tu))R_1 = q_{\w_1}(j^t(u))
$$
and the right hand side becomes
$$
(\varphi_{\w_2}(u) \otimes \varphi_{\w_2}(u))(F_{21}R_2F^{-1})=
q_{\w_2}(u)f(u,u)f(u,u)^{-1} = q_{\w_2}(u) \,.
$$
Therefore, $(\G^{\w_1},q_{\w_1})$ and $(\G^{\w_2},q_{\w_2})$ are equivalent.
Conversely, assume that there exists an isomorphism $j:\G^{\w_2} \map \G^{\w_1}$ such that
$q_{\w_1}j = q_{\w_2}$. Set
$$
\hat{j}(\infl {\w_1})(u,v,w) = \infl {\w_1}(j(u),j(v),j(w))
$$ and
$$
\hat{j}d_{\w_1}[u|v]=d_{\w_1}[j(u)|j(v)]
$$
for any $u,v \in \G^{\w_2}$. Then,
$(\hat{j}(\infl \w_1^{-1}), \hat{j}(d_{\w_1}))$ is also an abelian 3-cocycle of $\G^{\w_2}$.
The trace of this cocycle is $q_{\w_1} j$. By Theorem \ref{t3.2}, there exists a normalized
2-cochain $f$ on $\G^{\w_2}$ such that
\begin{eqnarray}
\infl \w_1^{-1}(j(u),j(v),j(w)) & =& \w_2^{-1}(u,v,w) \delta f(u,v,w) \label{eq8.3}\\
d_{\w_1}[j(u)|j(v)] & = & d_{\w_2}[u|v] f(u,v)/f(v,u)\label{eq8.4}
\end{eqnarray}
for any $u,v,w \in \G^{\w_2}$. Let
$$
F=\sum_{u,v \in \G^{\w_2}} f(u,v) E_u \otimes E_v \,.
$$
Notice that $\e(E_u) = \<E_u,\s(1,1)\> = \delta_{1,u}$ and $E_u E_v = \delta_{u,v}E_u$ for
any $u,v \in \G^{\w_2}$. This implies that $F$ is invertible in $\DG{2} \otimes \DG{2}$ and
$$
(id \otimes \e) F = 1_\DG{2}= (\e \otimes id) F\,.
$$
For simplicity, we denote the linear extension of $j$ from $\DG{2}_0$ to $\DG{1}_0$  by
the same symbol $j$. Obviously, $j$ is a bialgebra isomorphism. By Lemma \ref{l8.2},
$j^{t}$ is also a bialgebra isomorphism. For any $u,v,w \in \G^{\w_2}$,
\begin{eqnarray*}
(\varphi_{\w_2}(u) \otimes \varphi_{\w_2}(v) \otimes \varphi_{\w_2}(w))
(j^t \otimes j^t \otimes j^t) \Phi_1 & = &
(\varphi_{\w_1}(j(u)) \otimes \varphi_{\w_1}(j(v)) \otimes \varphi_{\w_1}(j(w)))\Phi_1 \\
& = & \infl {\w_1}^{-1} (j(u),j(v),j(w))
\end{eqnarray*}
and
$$
(\varphi_{\w_2}(u) \otimes \varphi_{\w_2}(v) \otimes \varphi_{\w_2}(w))
(F_{23} (id \otimes \Delta)(F)\Phi_2
(\Delta \otimes id)(F^{-1})F_{12}^{-1}) = (\infl {\w_2}) \delta f(u,v,w)\,.
$$
Hence, by equation (\ref{eq8.3}), we obtain
$$
(j^t \otimes j^t \otimes j^t)\Phi_1 =
F_{23} (id \otimes \Delta)(F)\Phi_2
(\Delta \otimes id)(F^{-1})F_{12}^{-1}) \,.
$$
It follows from equation (\ref{eq8.04}) that
$$
(\varphi_{\w_2}(u) \otimes \varphi_{\w_2}(v))(j^t \otimes j^t)R_1 =
(\varphi_{\w_1}(j(u)) \otimes \varphi_{\w_1}(j(v)))R_1= d_{\w_1}[j(v)|j(u)]
$$
and
$$
(\varphi_{\w_2}(u) \otimes \varphi_{\w_2}(v))F_{21}R_2F^{-1}= f(v,u)f^{-1}(u,v)
d_{\w_2}[v|u]\,.
$$
Hence, by equation (\ref{eq8.4}),
$$
(j^t \otimes j^t)R_1 = F_{21}R_2F^{-1}\,.
$$
Therefore, $j^t: \DG{1}\map \DG{2}_F$ is a quasi-bialgebra isomorphism. \qed

\begin{thm}\label{t8.5}
 Let $G_1$, $G_2$ be finite abelian groups and $\w_1 \in Z^3(G_1,\BC^*)_{ab}$,
$\w_2 \in Z^3(G_2,\BC^*)_{ab}$. Then, $D^{\w_1}(G_1)$ and $D^{\w_2}(G_2)$ are
gauge equivalent if, and only if
$\M{D^{\w_1}(G_1)}$ and $\M{D^{\w_2}(G_2)}$ are equivalent additive braided tensor
categories; that is, there exists a braided tensor equivalence
$({\cal F},\varphi_0,\varphi_2)$ between $\M{D^{\w_1}(G_1)}$ and $\M{D^{\w_2}(G_2)}$ with
$\cal F$ being additive.
\end{thm}
\pf The ``only if'' part of the statement is
well-known (cf. \cite{Kassel}). Conversely, let
$({\cal F},\varphi_0,\varphi_2)$ be an additive braided tensor equivalence from
$\M{K}$ to $\M{S}$ where $K=D^{\w_1}(G_1)$ and $S=D^{\w_2}(G_2)$. By Morita theory
(cf. \cite{Fuller}, \cite{Lam}),
${\cal F}$ induces an algebra isomorphism $\a$ of $\End_S({\cal F}(K))$ with $K$ and so ${\cal
F}(K)$ becomes a $S$-$K$-bimodule. Moreover, ${\cal F}$ is equivalent to the tensor functor
${\cal F}(K) \otimes_K \,?$. Note that both $S$ and $K$ are semisimple and their irreducible
modules are 1-dimensional. Therefore, ${\cal F}(K)= {_S\!S}$ and the right $K$-module action
is given by the isomorphism $\a^{-1}:K \map S$ where $S$ is naturally identified with
$\End({_S\! S})$. For any $M \in \M{K}$, define $\a^*(M)$ to be the left
$S$-module with the underlying space $M$ and the left $S$-action  given by
$$
s m = \a(s)m
$$
for any $s \in S$ and $m \in M$. It is easy to see that
$\Hom_K(M,N) = \Hom_S(\a^*(M),\a^*(N))$ for any $M,N \in \M{K}$. Hence, $\a^*$ defines an
additive functor from $\M{K}$ to $\M{S}$. It is straightforward to show that $\a^*$ and
${\cal F}(K)\otimes _K\, ?$ are equivalent functors. Without loss of generality, we may assume
${\cal F} = \a^*$.

Now $\BC \stackrel{\varphi_0}{\map} \a^*(\BC)$ is  scalar multiplication by a nonzero
complex number $\kappa$. For any $M \in \M{S}$, define $\eta_M$ to be  scalar
multiplication $\kappa^{-1}$. Then $\eta : \a^* \map \a^*$ is a natural isomorphism. Hence,
$(\a^*,id,\kappa \varphi_2)$ is a braided tensor equivalence from
$\M{K}$ to $\M{S}$. Consider $\kappa \varphi_2(K,K)$. Then
$\kappa\varphi_2(K,K)(1_K \otimes 1_K)$ is an invertible in $K \otimes K$ and we set
$F^{-1} = \kappa\varphi_2(K,K)(1_K \otimes 1_K)$. Then, $F$ is a gauge transform of $K$
(cf. \cite{Kassel},p381) and
\begin{eqnarray*}
\a(\Phi_S) &=&F_{23}(id \otimes \Delta_K)(F) \Phi_K (\Delta \otimes id)(F^{-1})F^{-1}_{12}\\
\a(R_S) & = & F_{21} R_K F^{-1}
\end{eqnarray*}
where $\Phi_S$, $R_S$ and $\Phi_K$, $R_K$ are the associators and $\cal R$-matrices of $S$ and
$K$ respectively. \qed

\begin{defn}
{\rm
Let $\G$ be a finite abelian group and $b:\G \times \G \map \BC^*$ a symmetric
bicharacter. For any subset $M$ of $G$,  denote by $M^\perp$ the subgroup
$$
\{x \in \G\,|\,b(x,M) = 1\}
$$
of $\G$. A subgroup $N$ of $\G$ is called a {\em metabolizer} of $(\G,b)$ if
$N=N^\perp$.
A quadratic map $q:\G \map \BC^*$ is called {\em non-degenerate} if the associated
bicharacter $b_q(x,y)=\frac{q(x y)}{q(x)q(y)}$ is non-degenerate. A subgroup $N$ of
$\G$ is called a {\em metabolizer} of $(\G,q)$ if $q|_N = 1$ and $N$ is a
metabolizer of $(\G,b_q)$. A metabolizer $N$ of $(\G,q)$ or $(\G,b)$ is called
{\em split} if
$N$ is a summand of $\G$.
}
\end{defn}
\begin{remark}\label{r8.6}\hfill
{\rm
\begin{enumerate}
\item[(i)] Let $G$ be a finite abelian group and $\w$ a normalized 3-cocycle in
$Z^3(G,\BC^*)_{ab}$. Then,
$q_\w(\widehat{G})=1$. Moreover, if $b_\w(\s(\a,x),\s(\l,1))=1$ for all $\l \in \widehat{G}$,
then $\l(x)=1$ for all $\l \in \widehat{G}$. This implies $x=1$. Therefore,
$\widehat{G}^\perp =\widehat{G}$ and hence $\widehat{G}$ is a metabolizer of $(\G^\w,q_\w)$.
\item[(ii)] If $(\G,q)$ is non-degenerate with metabolizer $G$, then $|G|^2=|\G|$.
\end{enumerate}
}
\end{remark}

Let $G$, $G'$ be finite abelian groups, $\w \in Z^3(G,\BC^*)_{ab}$ and $\w' \in
Z^3(G',\BC^*)_{ab}$. Let $\zeta \in Z^3(G_1 \times G_2,\BC^*)_{ab}$ be the product of the
inflations of $\w$ and $\w'$, namely
$$
\zeta =(\infl \w)(\infl \w')
$$
where the inflations are induced by the natural surjections $G\times G' \map
G$ and  $G\times G' \map
G'$. By the remarks following  Proposition \ref{pTP.1},
there is an isomorphism of groups $\iota : \G^{\w}(G) \times \G^{\w'}(G) \map \G^\zeta(G
\times G')$ given by
$$
\iota\left( \sum_{g \in G} \l_x(g) e(g) \otimes x, \sum_{h\in G'} \l'_y(h) e(h) \otimes
y \right) = \sum_{(g,h) \in G \times G'} \l_x(g)\l'_y(h) e(g,h) \otimes (x,y)
$$
where $\delta\l_x =\w_x$, $\delta\l'_y =\w'_y$ and $x \in G$, $y \in G'$. Since, for any
$u=\sum_{(g,h) \in G\times G'} \l_{(x,y)}(g,h) e(g,h) \otimes (x,y) \in \G^\zeta(G \times G')$,
$q_\zeta(u)=\l_{(x,y)}(x,y)$, $\iota$  defines an equivalence of the quadratic forms
$(\G^\zeta,q_{\zeta})$ and $(\G^\w ,q_\w) \perp (\G^{\w'} ,q_{\w'})$.

Conversely, suppose that
$H$ is a finite abelian group and $\eta \in Z^3(H,\BC^*)_{ab}$ such that $(\G^\eta,q_\eta)$ is
equivalent to the orthogonal sum  $(\G^\w,q_\w) \perp (\G^{\w'},q_{\w'})$. Then,
$(\G^\eta(H),q_\eta)$ and $(\G^\zeta(G\times G'),q_\zeta)$ are equivalent  quadratic
forms. By virtue of Proposition \ref{t8.4}, $D^\eta(H)$ is equivalent to $D^\w(H) \otimes
D^{\w'}(H')$ as quasi-triangular quasi-bialgebras. This proves

\begin{prop}\label{p8.7}
Let $G$ be a finite abelian group and $\eta \in Z^3(G,\BC^*)_{ab}$.  Then, $D^\eta(G)$ is
equivalent to
$D^\w(H) \otimes D^{\w'}(H')$ as quasi-triangular quasi-bialgebras for some abelian groups $H$,
$H'$ and
$\w \in Z^3(H,\BC^*)_{ab}$, $\w' \in Z^3(H',\BC^*)_{ab}$ if, and only if,
$(\G^\eta,q_\eta)$ is equivalent to $(\G^\w,q_\w) \perp (\G^{\w'},q_{\w'})$.
\end{prop}

%%%%%%%
%
%
%\input section9.tex
\section{Lattices}\label{s10}
We use the following notation: $M$ is a self-dual, even,  lattice
with respect to the nonsingular bilinear form
$$
\<\bullet,\bullet\> : M \times M \map \BZ\,.
$$
Thus $x \mapsto \<x,\bullet\>$ is an isomorphism of $M$ with
$\Hom_\BZ(M,\BZ)$, moreover $\<x,x\>$ is an even integer for $x \in M$. Note that
$\<\bullet,\bullet\>$ is not necessarily positive definite.

Let $E$ be the space $\BR \otimes_\BZ M$ equipped with the $\BR$-linear extension of
$\<\bullet,\bullet\>$, let $M \C L \C E$ with $|L :M| < \infty$, and let
$L_0 = \{x \in E\,|\,\<x,L\> \C \BZ\}$ be the $\BZ$-dual of $L$. We set $G =
L/M$. There is a short exact sequence
\begin{equation}\label{eq10.1}
0 \map M/L_0 \map L/L_0 \map L/M \map 0
\end{equation}
and we pick a section $s : L/M \map L$, such that
$s(0)=0$, which naturally defines a section
$\overline{s}: L/M \map L/L_0$.

Because $M$ is self-dual, the pairing $\<\bullet,\bullet\> : L \times L \map \BQ$
induces a perfect pairing
\begin{equation}\label{eq10.2}
p : M/L_0 \times L/M \map S^1\,,\quad (x+L_0,y+M) \mapsto e^{2\pi i\<x,y\>}
\end{equation}
and so there is a natural identification of $M/L_0$ with
$\Hom(L/M,S^1) = \widehat{G}$. Thus the sequence (\ref{eq10.1}) is of the type we
have been considering. That is, the triple $L_0 \C M \C L$ defines an element of
$H^2_{ab}(G,\widehat{G})$ with $G=L/M$, and it is well-known that a 2-cocycle
$\b \in Z^2_{ab}(G,\widehat{G})$ which corresponds to the triple is defined, using
the section $\overline{s}$, via
\begin{equation}\label{eq10.3}
\b(x,y)=\overline{s}(x)+\overline{s}(y)-\overline{s}(x+y)
\end{equation}
for $x$, $y \in G$.

Following Dong-Lepowsky \cite{DL}, we pick an alternating, bilinear map
$c: L \times L \map S^1$ with the property that
\begin{equation}\label{eq10.3.1}
c(x,y)=(-1)^{\<x,y\>}
\end{equation}
for $x,y \in M$. Such $c$ always exists (cf. \cite{DL}, Remark 12.18).
\begin{prop}\label{p10.1}
With the previous notation, set
\begin{equation}\label{eq10.4}
\w(g,x,y)= c(s(g),s(x)+s(y)-s(x+y))e^{\pi
i\<s(g),s(x)+s(y)-s(x+y)\>}
\end{equation}
\begin{equation}\label{eq10.5}
\t_g(x)= c(s(g),s(x))e^{\pi i\<s(g),s(x)\>}\,.
\end{equation}
Then $\w \in Z^3(G,\BC^*)_{ab}$, $\w_g =\delta \t_g$, and
$\Lambda_G([\w])=[\b]$.
\end{prop}
\pf After some computation, we find that for $g,h,x,y \in G$,
\begin{equation}\label{eq10.6}
\begin{split}
\delta \w (g,h,x,y)  = &
c(s(g)+s(h)-s(g+h),s(x)+s(y)-s(x+y)) \\
& \cdot e^{\pi i\<s(g)+s(h)-s(g+h),s(x)+s(y)-s(x+y)\>}\,.
\end{split}
\end{equation}
Note that both $a=s(g)+s(h)-s(g+h)$ and
$b=s(x)+s(y)-s(x+y)$ lie in
$M$. Using equation (\ref{eq10.3.1}), the expression of (\ref{eq10.6}) is equal to 1 and,
as a result, $\w$ is a normalized 3-cocycle.

The identity $\w_g  = \delta \t_g$ follows immediately from (\ref{eq10.4}) and
(\ref{eq10.5}). Finally, using the expression (\ref{eq2.8}) together with
(\ref{eq10.5}) shows that $\Lambda_G([\w])$ is the cohomology class represented
by $\b_1 \in Z^2_{ab}(G,\widehat{G})$ where
\begin{eqnarray*}
\b_1(x,y)(g)& = &
\frac{c(s(g),s(x))c(s(x),s(g))c(s(g),s(y))c(s(y),s(g))}
{c(s(g),s(x+y))c(s(x+y),s(g))}\, e^{2 \pi i
\<s(g),s(x)+s(y)-s(x+y)\>}
\end{eqnarray*}
for $x,y,g \in G$. Since $c$ is alternating, equations (\ref{eq10.2}) and
(\ref{eq10.3}) yields
$$
\b_1(x,y)(g) = p(\b(x,y),g)\,.
$$
So by definition of the pairing  $p$, $\b_1 = \b$. \qed\\

We now have two quadratic forms associated to this situation: $(L/L_0,q_L)$ where
$\mbox{$q_L(x+L_0)$}=e^{\pi i \<x,x\>}$ for $x \in L$; and $(\G^\w,q_\w)$ canonically associated
to the 3-cocycle $\w$ (cf. (\ref{eq8.1})).
\begin{lem}\label{l10.2}
These two quadratic spaces are equivalent.
\end{lem}
\pf Let notation be as before. From section \ref{s2}, $\G^\w$ is a central extension of $G$
by $\widehat{G}$. The additive version of the corresponding  2-cocycle  $\b$ is given by
(\ref{eq10.3}). Here, as before, $G=L/M$ and $\hat{G} = M/L_0$ under the identification
given in  (\ref{eq10.2}). There is an isomorphism $j : L/L_0 \map \G^\w$ defined for $x \in L$
via
\begin{equation}\label{eq10.7}
x + L_0 \mapsto \s(x+L_0-\overline{s}(x+M),x+M)\,.
\end{equation}
Now use equations (\ref{eq8.1}), (\ref{eq10.5}) to see that
\begin{equation}\label{eq10.8}
q_\w(j(x+L_0))=p(x+L_0-\overline{s}(x+M),x+M)e^{\pi i \<s(x),s(x)\>}\,.
\end{equation}
Write $s(x+M)=x+m$ for some $m \in M$. Then we see that
\begin{eqnarray*}
q_\w(j(x+L_0))& = & e^{\pi i (2 \<x-(x+m),x\> + \<x+m,x+m\>) }\\
              & = & e^{\pi i\<x,x\>} = q_L(x+L_0)\,.
\end{eqnarray*}
This completes the proof of the lemma. \qed\\

We now begin again, this time with a non-degenerate quadratic form $q$ on the finite abelian
group $\G$. It is convenient to write $q$ additively, that is $q:\G \map \BQ/\BZ$. We also
consider $q_1=2q$ as a quadratic form with values in $\BQ/2\BZ$. According to \cite{Wa63},
Theorem (6), there is a rational lattice $(H',\<\bullet,\bullet\>)$ i.e.
$\<\bullet, \bullet\> : H' \times H' \map \BQ$ non-degenerate, an even (integral) sublattice
$H \C H'$, and an equivalence of quadratic forms
$$
\xymatrix{\G \ar[r]^-{q_1} & \BQ/2\BZ \\
H'/H \ar[u]^\cong_j \ar[ru]_-{\overline{q}} & \\
}
$$
where $\overline{q}$ is the quadratic form arising from $\<\bullet,\bullet\>$. Thus, for
$x,y \in H'$ we have $\overline{q}(x+H) = \<x,x\> \pmod{2\BZ}$ and
\begin{eqnarray}
 q_1(j(x+H))& \equiv & \<x,x\> \pmod{2 \BZ}\\  \label{eq10.9}
b(j(x+H),j(y+H)) &\equiv & \<x,y\>  \pmod{ \BZ} \label{eq10.10}
\end{eqnarray}
where $b$ is the bilinear form on $\G$ associated to $q$.

Now assume that $(\G,q)$ has a metabolizer $\widehat{G}$, and write
$j^{-1}(\widehat{G}) =M/H$ for a sublattice $M$ with $H \C M  \C H'$. Since
$\widehat{G}^\perp = \widehat{G}$, it follows from (\ref{eq10.10}) that $M$ is a self-dual,
integral sublattice of $H'$. Similarly, equation (\ref{eq10.9}) shows that $M$ is an
even lattice. We now have
\begin{thm}\label{t10.1}
\begin{enumerate}
  \item[\rm (i)] Let $M \C L$ be a pair of rational lattices with
$M$ even and self-dual, and
  let $L_0 \C M$ be the $\BZ$-dual of $L$. Then the pair $(L/L_0,q_L)$ is a non-degenerate
  quadratic space with metabolizer $M/L_0$. Furthermore,
  \item[\rm (ii)] There is $\w \in Z^3(G,\BC^*)_{ab}$ such that the corresponding quadratic
  space $(\G^\w,q_\w)$ is equivalent to $(L/L_0,q_L)$. Moreover, $\Lambda(\w)$ is
cohomologous
  to the element of $Z^2_{ab}(L/M,M/L_0)$ defined by (\ref{eq10.3}).
  \item[\rm (iii)] Conversely, suppose that $(\G,q)$ is a non-degenerate quadratic space with
  metabolizer. Then there are lattices $M \C L$ as in (i) such that $(\G,q)$ and
  $(L/L_0,q_L)$ are equivalent quadratic spaces.
  \item[\rm (iv)] Every non-degenerate quadratic space $(\G,q)$ with a metabolizer
  $\widehat{G}$ is equivalent to one of the form $(\G^\w,q_\w)$ for suitable
  $\w \in Z^3(G,\BC^*)_{ab}$.
\end{enumerate}
\end{thm}
\pf (i) and (ii) follow from Proposition \ref{p10.1} and Lemma \ref{l10.2}, while (iii)
follows from the discussion following the proof of the Lemma. Part (iv) follows from
(i)--(iii). \qed

\begin{remark}\label{r10.1}
{\rm
>From Theorem \ref{t10.1} and our earlier work, we see that the following
pieces of
  data are more or less equivalent:
  \begin{enumerate}
    \item A pair of rational lattices $M \C L$ with $M$ even and self-dual, $L/M \cong G$.
    \item A non-degenerate quadratic space $(\G,q)$ with metabolizer $\widehat{G}$.
    \item A cohomology class $[\w] \in H^3(G,\BC^*)_{ab}$.
    \item A cohomology class $[\b] \in Z^2_{ab}(G,\widehat{G})$ which is in the image of $S$
    (cf. Theorem \ref{t6.4}).
  \end{enumerate}
}
\end{remark}
%%%%%%%%%%%%%%%
%
%
%\input section10.tex
\section{Gauss Sums}\label{s11}
We begin with a discussion of Gauss sums on pairs
$(M,q)$ consisting of a finite abelian group $M$
and a non-degenerate quadratic form $q:M\map
\BQ/\BZ$. Thus,
$q(x+y)=q(x)+q(y)+b(x,y)$ where $b: M\times M \map \BQ/\BZ$ is  a non-degenerate
bilinear form. In our applications, $(M,q)$ will be the pair $(\G^\w,q_\w)$ that we
have already considered, or something related to it. Our discussion is closely
related to that of tom Dieck (\cite{T}, 2.2) adapted to our present needs.

The Gauss sum $G(M,q)$ is defined via
$$
G(M,q)= \sum_{m \in M} e^{2\pi i q(m)}\,.
$$
The Gauss map\footnote{we say map because $\g$
can be considered as function on a suitable Witt
group.}
$\g$ is defined via
$$
\g(M,q) = \frac{1}{\sqrt{|M|}}G(M,q)\,.
$$
\begin{lem}\label{gs:l1}
Let $N \C M$ be a subgroup of $M$ such that the restriction of $q$ to $N$ vanishes
identically. Then $q$ induces a non-degenerate quadratic form $\overline{q}$ on
$N^\perp/N$, and we have
$$
\g(M,q)=\g(N^\perp/N,\overline{q})\,.
$$
\end{lem}
\pf Note that $\overline{q}$ is defined via
$\overline{q}(m+N)=q(m)$ for $m \in N^\perp$. It is easy to see that this is
well-defined and yields a quadratic form on $N^\perp /N$ which is non-degenerate.
Then we calculate:
\begin{eqnarray*}
\g(M,q)&=&
\frac{1}{\sqrt{|M|}}\frac{1}{|N|} \sum_{n \in N} \sum_{m \in M} e^{2 \pi i q(m+n)} \\
& = & \frac{1}{\sqrt{|M|}}\frac{1}{|N|} \sum_{m \in M} e^{2 \pi i q(m)}
\sum_{n \in N} e^{2\pi i b(m,n)}
\end{eqnarray*}
where we used $q(m+n)=q(m)+q(n)+b(m,n)$ and $q(n)=0$ for $n \in N$. Now the inner sum
vanishes whenever $m \not\in N^\perp$ by the orthogonality of group characters, and
is otherwise equal to $|N|$. So in fact
\begin{eqnarray*}
\g(M,q)&=&
\frac{1}{\sqrt{|M|}} \sum_{m \in N^\perp} e^{2 \pi i q(m)} \\
& = & \frac{|N|}{\sqrt{|N||N^\perp|}} \sum_{m \in N^\perp/N} e^{2 \pi i \overline{q}(m)}\\
& = & \g(N^\perp/N,\overline{q})\,. \quad\quad\quad\square
\end{eqnarray*}
\begin{cor}\label{gs:c1}
If $(M,q)$ has a metabolizer then $\g(M,q)=1$.
\end{cor}

\begin{lem}\label{gs:l2}
Let $p$ be a prime and $(M,q)$ a non-degenerate quadratic space on the non-trivial
cyclic $p$-group $M$. One of the following holds:
\begin{enumerate}
  \item[\rm (i)] $p$ is odd, $|M|$ is a square, and $(M,q)$ has a metabolizer.
  \item[\rm (ii)] $p\equiv 1 \pmod{4}$, $|M|$ is not a square, and $\g(M,q)=\pm1$\,.
  \item[\rm (iii)] $p \equiv 3 \pmod{4}$, $|M|$ is not a square, and $\g(M,q)= \pm i$\,.
  \item[\rm (iv)] $p=2$ and $\g(M,q)$ is a primitive 8th root of unity.
\end{enumerate}
\end{lem}
\pf Pick $N \C M$ maximal subject to the condition that the restriction of $q$ to $N$
is identically zero. Set $V= N^\perp/N$, so that $\g(M,q)=\g(V,\overline{q})$ as in
Lemma \ref{gs:l1}. Now the maximality of $N$ implies that $(V,\overline{q})$ is anisotropic,
 that is if $v \in V$ and $\overline{q}(v)=0$ then in fact $v=0$. It is then easy to see
 that $|V| \le p^2$, and even $|V| \le p$ if $p$ is odd. First assume that $p$ is
odd.
 Then (i) holds if $V=0$, and otherwise $|V|=p$ and
$|M|=|N|^2|V|$ is not a square. Then (ii) or (iii) holds by the classical calculation of
Gauss sums (cf. \cite{T}).

Now take $p=2$. If $|V|=2$ or 4 it is an easy direct calculation that (iv) holds. To complete
the lemma, we must show that, in fact, $V \ne 0$ if $p=2$. Assume, by way of contradiction,
that indeed $V=0$. Then $|M|=|N|^2=2^{2n}$, say. Let $x$ be a generator of $M$. Then $2^nx$
generates $N$ and therefore we have $q(2^nx)=0$.
As $b$ is non-degenerate, $b(x,x)=k/2^{2n}$ for some odd integer $k$. It follows from the
equality
$$
q(2x)-2q(x)=b(x,x)
$$
that $q(x)=k'/2^{n+1}$ for some odd integer $k'$. Hence $q(2^nx)=k'/2 \ne 0$, contradiction.
\qed\\

%Now we have $b(2^{2n}x,x)=0$ and
%$b(2^{2n-1}x,x)\ne 0$, the latter because $b(\bullet,\bullet)$ is non-degenerate. Hence we get
%$b(2^{2n-1}x,x)=1/2 = b(2^n x,2^{n-1}x)$ and then also $b(2^{n-1} x,2^{n-1}x)=\pm 1/4$. Now
%it follows that $0=q(2^{n-1} x+2^{n-1}x)= 2q(2^{n-1} x)+b(2^{n-1} x,2^{n-1}x)$, whence
%$q(2^{n-1} x)$ has order exactly 8.
%Since $2^{n} x+2^{n-1}x$ has the same order as $2^{n-1}x$ then $2^{n-1}x= k(2^n x +2^{n-1}x)$
%for some odd integer $k$. Thus $k^2 \equiv 1 \pmod{8}$ and therefore
%$q(2^{n-1} x)=q(2^{n} x+2^{n-1}x)$ since $q(2^{n-1}x)$ has order $8$.
%Finally, $q(2^nx+2^{n-1}x)=q(2^nx)+q(2^{n-1}x)+ b(2^{n} x,2^{n-1}x)$, and we conclude from
%prior assertions that $b(2^{n} x,2^{n-1}x)=0$, contradiction. \qed\\

Now consider a non-degenerate pair $(M,q)$ which has a metabolizer, and where $M$ is a $p$-group
with
$p \equiv 1 \pmod{4}$. We may write $M$ as an orthogonal direct sum of pairs
$(M_i,q_i)$ with $M_i$ cyclic (\cite{Wa63}, Theorem 4). By Lemma \ref{gs:l2}, if $|M_i|$
is a square  then $(M_i,q_i)$ has a metabolizer and $\g(M_i,q_i)=1$.
Since
$|M|$ is a square (by Remark \ref{r8.6} (ii)) and $\g(M,q)=1$, there must be a even number of
$M_i$ whose order is not a square and satisfy $\g(M_i,q_i)=-1$, and an even number of $M_i$
whose order is not a square and satisfy $\g(M_i,q_i)=1$. So we see that now we may write (with a
change of notation)
\begin{equation}\label{gs:eq1}
(M,q)=(M_1,q_1) \perp \cdots \perp (M_r,q_r)
\end{equation}
with $\g(M_i,q_i)=1$ for each $i$, and  $M_i$ is either cyclic of square order or the product
of two cyclic groups, each of non-square order.

Before we proceed further with our discussion, let us recall Wall's nomenclature for
non-degenerate symmetric bilinear forms on a finite cyclic group. Let $G$ be a cyclic group of
order $p^n$ with generator $x$, where $p$ is an odd prime. A bilinear form $b\:G \times G \map
\BQ/\BZ$ with
$b(x,x)=\e p^{-n}$ is of {\em type A} if $\e=1$, and it is of {\em type B} if $\e$ is a quadratic
non-residue modulo $p$ (cf. \cite{Wa63}).

We assert that each $(M_i,q_i)$ in (\ref{gs:eq1}) has a cyclic
metabolizer. This has already been explained if
$M_i$ is cyclic. Otherwise, $M_i = H \oplus K$ with $H$, $K$ both non-zero cyclic groups
and $\g(H,q_i)\g(K,q_i)=1$. Then $\g(H,q_i)=\g(K,q_i)$ and from section 5 of \cite{Wa63} we
see that $H$ and $K$ are either both of type A or both of type B. More
precisely, there are generators $x,y$ of $H,K$ respectively (and of order $p^h$, $p^k$
respectively, say) such that $b_i(x,x)=\e p^{-h}$ and $b_i(y,y)=-\e p^{-k}$ where $\e$ is either
1 or  a non-residue $\pmod{p}$. (Here we have used the fact that
$\left( \frac{-1}{p}\right)=1$ for $p \equiv 1\pmod{4}$.) Now if $h=r-u$ and $k=r+u$ for
nonnegative integers $r,u$ we find that the subgroup generated by $x+p^uy$ is the metabolizer
we require.

We can rewrite (\ref{gs:eq1}) in terms of twisted quantum doubles using Proposition \ref{p8.7}.
As each $(M_i,q_i)$ has a cyclic metabolizer $G_i$, say, then we know (cf. Remark
\ref{r10.1}) that the pair
$(M_i,q_i)$ can be realized via $D^{\w_i}(G_i)$ for suitable $\w_i \in Z^3(G,\BC^*)_{ab}$.
That is, $(M_i,q_i)$ is equivalent to $(\G^{\w_i},q_{\w_i})$.
The orthogonal sum (\ref{gs:eq1}) corresponds to the tensor product of twisted doubles. So
if we start with a pair $(G,\w)$ which gives rise to $(M,q)$, we can conclude

\begin{thm}\label{gs:t1}
Suppose that $p \equiv 1 \pmod{4}$ is a prime, $G$ an abelian $p$-group, and
$\w \in Z^3(G,\BC)_{ab}$. Then there are cyclic $p$-groups $G_i$ and cocycles
$\w_i \in Z^3(G_i,\BC)_{ab}$, $1 \le i \le r$, such that the two quasi-triangular,
quasi-Hopf algebras $D^\w(G)$ and $\bigotimes\limits_{i=1}^r D^{\w_i}(G_i)$ are
gauge equivalent. \qed
\end{thm}

\begin{remark}
{\rm
It is evident from our proof that the choice of the isomorphism types of the $G_i$ is
far from unique. For example, suppose that $M$ is the orthogonal sum of four cyclic group of
orders $p$, $p^3$, $p^5$, $p^7$, and each of type A. Then our proof shows that $D^\w(G)$ is
equivalent to $D^{\w_1}(G_1) \otimes D^{\w_2}(G_2)$, as quasi-triangular quasi-bialgebras, where
$G_1$,
$G_2$ can be chosen to be have orders $p^2$, $p^6$; $p^3$, $p^5$;  or $p^4$, $p^4$. Of course,
the number
$r$ of tensor factors that occur depends only on $(G,\w)$, and we always have
$|G|= \prod\limits_{i=1}^r |G_i|$.
}
\end{remark}

Theorem \ref{gs:t1} fails for primes not congruent to $1\pmod{4}$. To treat this case we
need

\begin{lem}\label{gs:l3}
Let $(M,q)$ be a non-degenerate quadratic space on the $p$-group $M$. Suppose that $|M|$ is a
square. Then $\g(M,q)=1$ if, and only if, $(M,q)$ has a metabolizer.
\end{lem}
\pf
Let $N$ and $V$ be as in the proof of Lemma \ref{gs:l2}. After Corollary \ref{gs:c1} we may
assume $\g(M,q)=1$ and $V\ \ne 0$ and try to reach a contradiction.

We see that $V$ is not cyclic by Lemma \ref{gs:l2}, and since $V$ is anisotropic then
$V \cong \BZ_p \oplus \BZ_p$, as is easily verified. From Lemma \ref{gs:l2} (and its proof),
we see that $\g(M,q)=-1$ if either $p$ is odd, or if $p=2$ and $V$ is not an orthogonal
sum of two cyclic groups. Otherwise, $p=2$, $V=\<x,y \>$, $b(x,y)=0$, $q(x)=q(y)=\pm 1/4$,
$q(x+y)=1/2$. Hence $\g(M,q)= \pm i$ in this case. \qed\\

Now consider the analogue of (\ref{gs:eq1}) in the case $p \equiv 3 \pmod{4}$. The arguments
preceding (\ref{gs:eq1}) together with Lemma \ref{gs:l2} show that we may still write
$(M,q)$ in the form of equation (\ref{gs:eq1}), with each $\g(M_i,q_i)=1$, and $|M_i|$ a
square. However, $M_i$ may now be a product of 1, 2 or 4 cyclic subgroups. By Lemma \ref{gs:l3},
each $(M_i,q_i)$ has a metabolizer. Thus there is an analogue of Theorem \ref{gs:t1} for
$p \equiv 3 \pmod{4}$, except that the groups $G_i$ cannot necessarily be chosen to be
cyclic. Rather, we have the property that the group $\G^{\w_i}(G_i)$ is a product of 1, 2
or 4 cyclic groups.

Finally, take $p=2$, and assume again that $\g(M,q)=1$.  We can decompose $(M,q)$ as follows
(\cite{Wa63}, section 5,6)
\begin{equation}\label{gs:eq2}
(M,q)=(P_1,q_1) \perp \cdots \perp (P_t,q_t)
\end{equation}
with each $P_i$ homogeneous and the product of at most 2 cyclic groups. From Lemma
\ref{gs:l2}, and the proof of Lemma \ref{gs:l3}, we see that the number of cyclic
factors in (\ref{gs:eq2}) is necessarily even. Hence $M$ is the product of
an even number of cyclic groups, which thus gives a different proof of Theorem
\ref{t6.5}. The analogue of Theorem \ref{gs:t1} for $p=2$ is the following: $D^\w(G)$ is gauge
equivalent to $\bigotimes\limits_{i=1}^r D^{\w_i}(G_i)$ and for each $i$, $\G^{\w_i}(G_i)$ is
the product of 1,2,4 or 8 cyclic groups.
%%%%%%%%%
%
%
%\input \section11.tex
\section{Metabolic Forms And Homogeneous Groups}\label{s12}
A {\em metabolic
form} is a triple $(\G,b,G)$ with $b: \G \times \G \map \BQ/\BZ$ a non-degenerate symmetric
bilinear form and $G$ a metabolizer of $(\G,b)$.
If $\G' \le  G$ and $G' \le G$ such that $(\G',b,G')$ also a metabolic form, then
$(\G',b,G')$ is called a sub-metabolic form of $(\G,b,G)$. Obviously, $(\{0\},b,\{0\})$ is
a trivial sub-metabolic form of $(\G,b,G)$. A metabolic form is called
{\em simple} if $(\G,b,G)$ is the only nontrivial sub-metabolic form.

A group homomorphism  $j\:\G_1\map\G_2$
is called a {\em morphism} mapping $(\G_1,b_1,G_1)$ to $(\G_2,b_2,G_2)$ if
$j(G_1) \C G_2$ and
$$
b_2(j(u),j(v))= b_1(u,v)
$$
for any $u,v \in \G_1$. If the morphism $j$ is bijective, then we call the two
metabolic forms  {\em equivalent} and denote this by
$$
(\G_1,b_1,G_1)\cong (\G_2,b_2,G_2)\,.
$$

\begin{prop}\label{p11.09}
Let $G$, $G'$ be finite abelian groups of odd order, $\w \in Z^3(G,\BC^*)_{ab}$
and $\w' \in Z^3(G',\BC^*)_{ab}$. Then $(\G^{\w}(G),b_{\w},\widehat{G}) \cong
(\G^{\w'}(G'),b_{\w'},\widehat{G}')$
if, and only if, $\w'$ and $\widehat{\rho} w$ are
cohomologous for some isomorphism $\rho :G \map G'$ where $\widehat{\rho}\w$ is as in Remark
\ref{r0.1}(iii).
%\begin{equation}\label{eq11.09}
%\widehat{\rho} \w (x,y,z)=\w (\rho^{-1} x, \rho^{-1} y, \rho^{-1} z)
%\end{equation}
%for $x,y,z \in G'$.
\end{prop}
\pf Let $\rho : G \map G'$ be an isomorphism such that
$\widehat{\rho}\w=\w' \delta b'$ for some normalized 2-cochain $b'$ on $G'$.
Let $\t \in T(\w)$. Then $\t'\in T(\w')$ where
$$
\t'_{x'}(y') = \t_{\rho^{-1} x'}(\rho^{-1} y') \frac{b'(y',x')}{b'(x',y')}
$$
for any $x',y, \in G'$.
The map $j :\H \map D^{\w'}(G')$ defined by
$$
j: e(g) \otimes x \mapsto  \frac{b'(\rho g, \rho x)}{b'(\rho x, \rho g)}\,
e(\rho g) \otimes \rho x
$$
is a bialgebra isomorphism. In particular, $j(\G^\w(G)) = \G^{\w'}(G')$. Moreover,
for any $\s(\a,x) \in \G^\w(G)$,
$$
j(\s(\a,x))= \sum_{g\in G}\a(g)\t_x(g) \frac{b'(\rho g, \rho x)}{b'(\rho x, \rho g)}\,
e(\rho g) \otimes \rho x
$$
and so
$$
q_{\w'}(j(\s(\a,x))= \a(x) \t_x(x) =q_\w(\s(\a,x))\,.
$$
It is obvious that $j(\widehat{G})=\widehat{G'}$ and hence $j$ is an
isomorphism of $(\G^{\w}(G),b_{\w},\widehat{G})$ with
$(\G^{\w'}(G'),b_{\w'},\widehat{G}')$.

Conversely, assume that $j$ is an isomorphism of
$(\G^{\w}(G),b_{\w},\widehat{G})$ with
$(\G^{\w'}(G'),b_{\w'},\widehat{G}')$. For $x \in G$,
$$
j\s_{\t}(1,x) = \s_{\t'}(b_x, \rho x)
$$
for some $b_x \in \widehat{G}'$ and $\rho x \in G$.
As $j(\widehat{G}) = \widehat{G}'$, the map $\rho: x \mapsto \rho
x$  defines an isomorphism from $G$ to $G'$. As
$$
b_{\w'}(\s_{\t'}(b_x, \rho x),\s_{\t'}(b_y, \rho y)) =  b_\w(\s_{\t}(1,x), \s_{\t}(1,y)),
$$
we have the equality
\begin{equation}\label{eq11.091}
b_x(y)b_y(x) \t'_{\rho x}(\rho y)\t'_{\rho y}(\rho x) = \t_x(y)\t_y(x)\,.
\end{equation}
Hence,
\begin{eqnarray*}
\frac{\t_x(z)\t_y(z)}{\t_{xy}(z)}{\w_z(x,y)} & = &
\frac{\t_x(z)\t_z(x)\t_y(z)\t_z(y)}{\t_{xy}(z)\t_z(xy)} \\
 &=& \frac{\t'_{\rho x}(\rho z)\t'_{\rho z}(\rho x) \t'_{\rho y}(\rho z)
\t'_{\rho z}(\rho y)}{\t'_{\rho (xy)}(\rho z)\t'_{\rho z}(\rho (xy))}
\frac{b_x(z)b_y(z)}{b_{xy}(z)}\\
 &=& \frac{\t'_{\rho x}(\rho z)\t'_{\rho y}(\rho z)}
{\t'_{\rho (xy)}(\rho z)}\w'_{\rho z}(\rho x,\rho y) (\delta b)(x,y)(z)
\end{eqnarray*}
where $b \in C^1(G,\widehat{G})$ given by $b(x,y)=b_x(y)$.
By equation (\ref{eq11.091}) and (\ref{eq2.8}),  we see that
$\Lambda(\widehat{\rho}\w)$ and $\Lambda(\w')$ are cohomologous.
As $|G|$ is odd, it follows from Corollary \ref{c4.4} that
$[\widehat{\rho}\w] = [\w']$. \qed\\

For  two metabolic forms $(\G_1,b_1,G_1)$ and $(\G_2,b_2,G_2)$, one can define the
{\em orthogonal sum}, $(\G_1,b_1,G_1) \perp (\G_2,b_2,G_2)$ to be the metabolic form
$(\G_1 \times G_2, b, G_1\times G_2)$ where
$$
b((u_1,v_1),(u_2,v_2)) = b_1(u_1,v_1)+b_2(u_2,v_2)
$$
for any $u_1,u_2 \in G_1$ and $v_1,v_2 \in G_2$.

\begin{prop}\label{p11.1}
Let  $(\G,b,G)$ be a metabolic form and $(K,b|_K,H)$ a sub-metabolic form of
$(\G,b,G)$. Then, $(K^\perp, b|_{K^\perp}, G\cap K^\perp)$ is also a
sub-metabolic form of $(\G,b,G)$. Moreover,
$$
(\G,b,G)  \cong (K,b|_K,H) \perp (K^\perp, b|_{K^\perp}, G\cap K^\perp)\,.
$$
\end{prop}
\pf
As $b$ is non-degenerate on $K$, there exists a subgroup $Q < K^\perp$
such that $K +G = K +Q$.
So for any $x \in G$, $x = u+ v$ for some $u \in K$ and $v \in Q$. Now $H \C G$ and
$b(v,H)=0$. This implies $b(u,H)=0$. Since $H$ is a metabolizer of $(K,b|_{K})$,
$u \in H$, it follows  that $v \in G$ and hence $G  = H +(K^\perp \cap G)$.
Moreover,
$$
 (G \cap K^\perp)^\perp\cap K^\perp = (G+ K) \cap K^\perp=  (K+(K^\perp \cap G)) \cap
K^\perp = G \cap K^\perp \,.
$$
Hence $G \cap K^\perp$ is a metabolizer of $(K^\perp, b|_{K^\perp})$. It is straightforward to
show that
$$
(K,b|_K,H) \cong (\G,b|_{\G},G) \perp (K^\perp ,b|_{K^\perp},G \cap K^\perp)\cong
(\G,b,G)
$$
with the morphism $j:K \times K^\perp \map \G$ given by $j(u,v)=u+v$ for any
$(u,v) \in K \times K^\perp$. \qed

\begin{prop}\label{p11.2}
Every metabolic form is equivalent to an orthogonal sum of simple metabolic forms.
\end{prop}
\pf The result follows directly from Proposition \ref{p11.1} by induction. \qed

\begin{remark}
{\rm
It follows easily from Proposition \ref{p11.1} that a metabolic form $(\G,b,G)$ is simple
if $G$ is cyclic.
}
\end{remark}

We will denote by $(K,b|_K,H)^\perp$ the metabolic form
$(K^\perp,b|_{K^\perp},K^\perp \cap G)$ for any sub-metabolic form $(K,b|_K,H)$ of
$(\G,b,G)$.

\begin{lem}\label{l11.3}
Let $(\G,b,G)$ be a metabolic form with $G$ homogeneous of
exponent $p^n$. For any non-degenerate cyclic subgroup $D$ of $\G$
of order greater than $p^n$,
there exists a cyclic subgroup $H\C G$
such that $(H+D,b|_{H+D}, H)$ is a sub-metabolic form
of $(\G,b,G)$. Moreover, for any generator $y$ of $D$, there is a generator $x$ of $H$
satisfying
$$
x = x'+ {p^a}y
$$
for some $x' \in D^\perp$ where $|D|=p^{n+a}$.
\end{lem}
\pf Let $D$ be as in the statement of the lemma. Set
$G_0 = G \cap D^\perp$. Then $G/G_0$ is cyclic since the map
$g+G_0 \mapsto b(g,\bullet)$ defines an embedding from $G/G_0$ into $\widehat{D}$.
If $|G:G_0| <p^n$,  then $\Omega_p(G) \C G_0$. Therefore, $\Omega_p(G) \cap D$ is trivial
and so is $G \cap D$. But this implies $|D| \le p^n$ since
$$
 D/(G \cap D) \cong (G+D)/G \C \G/G \cong \widehat{G}.
$$
Hence $|G:G_0| = p^n$. Then, $G=G_0\oplus H$ for some
$H \le G$  is isomorphic to $\BZ_{p^n}$. $H$ can be chosen to
contain $G \cap D$ as  $G$ is homogeneous. In particular, we have
\begin{equation}\label{eq11.1}
G \cap D = H\cap D\,.
\end{equation}
Set $K=H+D$. Then $G_0 \C K^\perp$. Therefore,
$$
K^\perp \cap K = (K +K^\perp)^\perp \C (D+ H +G_0)^\perp = G^\perp \cap D^\perp =G_0 \,.
$$
As $G \cap K = H+(G \cap D) = H$ by equation (\ref{eq11.1}), $G_0 \cap K$ is trivial.
Hence, $K^\perp
\cap K=\{0\}$, i.e. $(K,b|_K)$ is non-degenerate. Moreover,
$$
H^\perp \cap K = H +(H^\perp \cap D) = H +((G_0+H)^\perp \cap D) = H+(G \cap D),.
$$
By (\ref{eq11.1}), $H^\perp \cap K = H$ and hence $H$ is a metabolizer of $(K,b|_K)$.

Let $\<x\>=H$ and $\<y\>=D$. Then, $x=x'+{s p^a}y$ for some integer $s$ and
$x' \in D^\perp$. If $p|s$, $b(p^{n-1}x,y)= 0$ and so $p^{n-1}x \in G_0 \cap H$ which
contradicts  $\ord(x) =p^n$. So there is $s' \in \BZ$ such that $ss'\equiv 1 \pmod{p}$. The
result follows if one replace $x$ by $s'x$. \qed
\begin{lem}\label{l11.4}
Let $p$ be an odd prime and $(\G,b,G_1)$,  $(\G,b,G_2)$  metabolic forms with
$G_1$ and $G_2$ isomorphic to $(\BZ_{p^n})^k$ for some positive integer $k$. Then, there exist
non-trivial sub-metabolic forms $(K_i,b|_{K_i}, H_i)$ of $(\G,b,G_i)$, $i=1,2$ such that
$$
(K_1,b|_{K_1}, H_1) \cong (K_2,b|_{K_2}, H_2)\,.
$$
\end{lem}
\pf Pick $D \le \G$ cyclic of maximal order such that $b|_D$ is non-degenerate
(cf. \cite{Wa63}). Then $|D| \ge p^n$. If $|D|=p^n$, $G_1$ and $G_2$ are  split metabolizers of
$(\G,b)$. It is straightforward to show that there exist $x_i \in G_i$ and $u_i, v_i \in \G$
such that $b(u_i,u_i)=1/p^n$, $b(v_i,v_i)=-1/p^n$, $b(u_i,v_i)=0$ and $x_i=u_i+v_i$
(cf. \cite{HM}). In particular, $(K_i,b|_{K_i}, \<x_i\>)$ is a sub-metabolic form of
$(\G,b,G_i)$ where $K_i$ is the subgroup generated by $u_i,v_i$. Moreover,
the map $K_1 \map K_2$, $u_1 \mapsto u_2$ and $v_1 \mapsto v_2$ defines an isomorphism of
$(K_1,b|_{K_1}, \<x_1\>)$ with $(K_2,b|_{K_2}, \<x_2\>)$.

Let $y$ be a generator of $D$ and $\ord(y)=p^{n+a}$ for some positive integer $a$. By Lemma
\ref{l11.3}, for $i=1,2$ there exists $H_i \le G_i$ such that
$(K_i,b|_{K_i},H_i)$ is a sub-metabolic form of $(\G,b,G_i)$ where $K_i=H_i+D$.
Moreover $H_i$ admits a generator $x_i$ satisfying
$$
x_i = x_i' + p^a y
$$
for some $x_i' \in D^\perp$. The map $K_1 \map K_2$, given by
$y \mapsto y$ and $x_1' \mapsto x_2'$, defines an isomorphism of
$(K_1,b|_{K_1},H_1)$ with $(K_2,b|_{K_2},H_2)$. \qed

\begin{thm}\label{t11.5}
Let $p$ be an odd prime and $(\G,b,G)$ a metabolic form with $G \cong (\BZ_{p^n})^k$.
\begin{enumerate}
  \item[\rm (i)]  $(\G,b,G)$ is equivalent to an orthogonal sum of sub-metabolic forms
with  metabolizers isomorphic to $\BZ_{p^n}$.
  \item[\rm (ii)] Let $(\G',b',G')$ be a metabolic form such that $G'\cong G$.
   Then $(\G',b',G') \cong (\G,b,G)$ if,
and only if $(\G',b')$ and $(\G,b)$ are equivalent.
\end{enumerate}
\end{thm}
\pf Statement (i) follows easily from Lemma \ref{l11.4} and Proposition \ref{p11.1} by
induction on $k$. For statement (ii), it is obvious that the equivalence of two metabolic
form implies the equivalence of the underlying bilinear forms. Conversely, we proceed by
induction on $k$. Let $j$ be an
isomorphism of $(\G',b',G')$ with $(\G,b,G)$. Then, $(\G,b,j(G'))$ is a metabolic form. By
Lemma \ref{l11.4}, there exist a nontrivial sub-metabolic form
$(K',b'|_{K'},H')$ of $(\G',b',G')$ and  a sub-metabolic form $(K,b|_{K},H)$ of
$(\G,b,G')$ such that
$$
(K',b|_{K'},H')\cong (K,b|_{K},H)\,.
$$
Thus, $({K'}^\perp,b'|_{{K'}^\perp})$ and $(K^\perp,b|_{K^\perp})$ are equivalent bilinear form.
By induction assumption $(K',b|_{K'},H')^\perp \cong (K,b|_{K},H)^\perp$ and hence
$$
(\G',b',G,) \cong (K',b|_{K'},H') \perp (K',b|_{K'},H')^\perp \cong
(K,b|_{K},H)\perp (K,b|_{K},H)^\perp \cong (\G,b,G)\,.
$$
\hfill\qed
\begin{thm}\label{t11.6}
Let $p$ be an odd prime with $G$  isomorphic to $(\BZ_{p^n})^k$.
\begin{enumerate}
  \item[\rm (i)] For any $\w \in Z^3(G,\BC^*)_{ab}$, there exists
  $\eta_1,\dots ,\eta_k  \in Z^3(\BZ_{p^n},\BC^*)_{ab}$ such that
  $D^\w(G)$ is gauge equivalent to $\bigotimes\limits_{i=1}^k D^{\eta_i}(\BZ_{p^n})$ as
  quasi-triangular quasi-bialgebras.
  \item[\rm (ii)] For any $\w,\w' \in Z^3(G,\BC^*)_{ab}$,
  $\H$ and $D^{\w'}(G)$ are gauge equivalent as quasi-triangular quasi-bialgebras if, and only
if, there exists
  $\rho \in \Aut(G)$ such that $[\widehat{\rho}\w]=[\w']$.
\end{enumerate}
\end{thm}
\pf (i) Consider a metabolic form  $(\G^\w,b_\w,\widehat{G})$ associated to $\w$.
By Theorem \ref{t11.5}, $(\G^\w,b_\w,\widehat{G})$ equivalent to the orthogonal sum
\begin{equation}\label{eq11.2}
(\G_1,b_1,\BZ_{p^n}) \perp \cdots \perp (\G_k,b_k,\BZ_{p^n})\,.
\end{equation}
By  Theorem \ref{t10.1} (iv), there exist $\eta_i \in Z^3(\BZ_{p^n},\BC^*)$ such that
$(\G_i, b_i)$ is equivalent to $(\G^{\eta_i},b_{\eta_i})$. Hence, by equation (\ref{eq11.2}),
the result follows. \\
(ii)  If $D^\w(G)$ and $D^{\w'}(G)$ are equivalent as quasi-triangular quasi-bialgebras, by
Theorem \ref{t8.4}, $(\G^\w,b_\w)$ and $(\G^{\w'}, b_{\w'})$ are equivalent bilinear forms.
It follows from Theorem \ref{t11.5}(ii) that $(\G^{\w}, b_{\w}, \widehat{G})$
and $(\G^{\w'}, b_{\w'}, \widehat{G})$ are equivalent metabolic forms. By
 Proposition \ref{p11.09}, $[\widehat{\rho}\w]=[\w']$ for some $\rho \in \Aut(G)$.
The ``only if'' part follows from Remark \ref{r0.1} (iii). \qed
%%%%%%%%%%%
%
%
%\input section12.tex
\section{Duality and Symmetry}\label{s13}
In this section we discuss various kinds of relations that exist between suitably chosen
$\H$ and $D^{\w'}(G')$. Our discussion is meant to illustrate the possibilities, and is
by no means exhaustive. We make use of Wall's classification of symmetric bilinear forms
\cite{Wa63} as well as results previously established in the present paper.

The first type relation we discuss--{\em symmetry}--arises from the possibility that
a quadratic space $(\G,q)$ may have metabolizers $G$, $G'$ which are not conjugate in the
corresponding orthogonal group. Indeed, they may not even be isomorphic. This is, in fact,
a rather common phenomena.

As an example, take $\G \cong \BZ_{p^2} \times \BZ_{p^2}$ ( and for convenience, $p$ an
odd prime). By \cite{Wa63}, the possible non-degenerate bilinear forms on $\G$ can be taken
to be the following: $A_{p^2} \oplus A_{p^2}$, $A_{p^2} \oplus B_{p^2}$. Of these, only one
has a (necessarily split) cyclic metabolizer:  the first if $p \equiv 1 \pmod{4}$ and the second
if
$p \equiv 3 \pmod{4}$. On the other hand, the subgroup of $\G$ consisting of elements of order
at most $p$ is a metabolizer in all cases. As we know, each of these metabolizers $G$
determines an abelian 3-cocycle $\w \in Z^3(G,\BC^*)_{ab}$ such that $\G^\w \cong \G$, and
$\w$ is even a coboundary if $G$ is cyclic. More is true: if they determine the same pair
$(\G,q)$ then the corresponding module categories are braided monoidally equivalent. We
thus conclude:

\begin{example}\label{ex12.1}
{\rm
Let $p$ be an odd prime, $G \cong \BZ_{p^2}$, $G' \cong \BZ_{p} \times \BZ_{p}$. Then there
is an abelian 3-cocycle $\w' \in Z^3(G',\BC^*)_{ab}$ such that $\M{D(G)}$ and
$\M{D^{\w'}(G')}$ are equivalent as braided tensor categories. The same argument applies if we
take $G' \cong \BZ_{p^n} \times \BZ_{p^n}$ and $G \cong \BZ_{p^{2n}}$, and to a host of
other situations.
}
\end{example}

Next we discuss a {\em duality} between the module categories corresponding to the twisted
double of homogeneous $p$-groups.

\begin{thm}\label{t12.2}
Let $n,k$ be a pair of positive integers, and let $G\cong \left(\BZ_{p^n}\right)^k$ and
$G' \cong  \left(\BZ_{p^k}\right)^n$ ($p$ an odd prime). Then the following hold:
\begin{enumerate}
\item[\rm (i)] There are exactly ${n+k \choose k}$ equivalence classes of  monoidal
categories of
 the form $\M{\H}$, for some $\w \in Z^3(G,\BC^*){a b}$.
\item[\rm (ii)] There is a canonical bijection between equivalence classes of monoidal
categories of the form $\M{\H}$ and those of the form $\M{D^{\w'}(G')}$.
\item[\rm (iii)] There are natural bijections between equivalence classes of braided monoidal categories of
the form $\M{\H}$ and those of the form $\M{D^{\w'}(G')}$.
\end{enumerate}
\end{thm}

\begin{example}\label{ex12.2}
{\rm
One can also count the number of braided monoidal categories (up to equivalence) of the form
$\M{\H}$, where $G=\left(\BZ_{p^n}\right)^k$ and $\w \in Z^3(G,\BC^*)_{ab}$. Thus we have
\begin{center}
\begin{tabular}{c||c}
$k$ & \# braided monoidal categories\\ \hline
$1$ & $(2n+1)$ \\
$2$ &  $2n^2+2n+1$\\
$3$ & $1/3 (4n^3+6n^2+8n+3)$ \\
$\cdots$ & $\cdots$
\end{tabular}
\end{center}
}
\end{example}

To prove these assertions, recall from Theorem \ref{t11.5} that if
$G = \left(\BZ_{p^n}\right)^k$ is homogeneous, $p$ odd, then the pair $(G,\w)$,
$\w\in Z^3(G,\BC^*)_{ab}$, determines a non-degenerate bilinear form $(\G,b)$ of the
form
\begin{equation}\label{eq12.1}
(\G,b) = (H_1,b_1) \perp \cdots\perp (H_k,b_k)
\end{equation}
where $(H_i,b_i)$ has a cyclic metabolizer $G_i$, $G= G_1 \times \cdots \times G_k$,
and $H_i$ is isomorphic to  $\BZ_{p^{n+a}} \bigoplus \BZ_{p^{n-a}}$ for some integer $a$
with $0\le a \le n$. Of course, $\G$ is isomorphic to the group of fusion rules $\G^\w$.
Denote by $\G_a$ the group $\BZ_{p^{n+a}} \bigoplus \BZ_{p^{n-a}}$ (where now $n$ is
fixed). If $\G_a$ occurs with multiplicity
$m_a$ in (\ref{eq12.1}), we can represent (\ref{eq12.1}) symbolically in the form
\begin{equation}\label{eq12.2}
\G = \perp_{a=0}^{n} m_a \G_a
\end{equation}
at least as far as the fusion rules are concerned. It is easy to see that we have
\begin{equation}\label{eq12.3}
\sum_{a=0}^n m_a =k\,.
\end{equation}
Furthermore, any $n+1$-tuple $(m_0, m_1,  \dots, m_n)$ of nonnegative integers satisfying
(\ref{eq12.3}) may be chosen in (\ref{eq12.2}) and which correspond to a pair $(\G,b)$ with
metabolizer $G$.
Thus the number of non-isomorphic groups of fusion rules $\G$ that correspond to $G$ is the
number of non-negative $(n+1)$-tuples $(m_0, m_1, \dots, m_n)$ satisfying (\ref{eq12.3}).
This is easily seen to equal ${n+k \choose k}$. Now part (i) of the Theorem follows from
Theorem \ref{t7.3}.

Note that if also $G' =  \left(\BZ_{p^k}\right)^n$ then by part (i), there are
 ${n+k \choose k}$ equivalence classes of tensor categories both of the type
 $\M{\H}$ and $\M{D^{\w'}(G')}$. We will establish a canonical bijection between these
two sets.

With $(m_0,m_1 ,\dots, m_n)$ as above, let $\l$ be the partition which contains the integer
$i$ with multiplicity $m_i$, $1 \le i \le n$. Thus $\l$ is a partition of the integer
$r = \sum_{i=1}^n im_i$. Note that $r \le nk$ by (\ref{eq12.3}). Let $\l^t$ be the dual
(conjugate) partition of $\l$, and let $s_j$ be the multiplicity of $j$ in $\l^t$ for
$j \ge 1$. So also $ r = \sum_{j=1}^k j s_j$, the upper limit $k$ in the sum arising
from the fact that the maximal part of $\l^t$ (i.e., the largest $j$ with $s_i > 0$) is
equal to the number of non-zero multiplicities $m_i$, and this is at most $k$ by
(\ref{eq12.3}).

Now observe that $\sum_{j=1}^k s_j =$ \# of parts of the partition
$\l^t= \max\{i|m_i > 0\} \le n$. Thus if we define $s_0 = n-\sum_{j=1}^ks_j$ then $s_j \ge 0$
for $j \ge 0$ and
\begin{equation}\label{eq12.4}
\sum_{j=0}^k s_j =n \,.
\end{equation}
Consider
\begin{equation}\label{eq12.5}
\G' = \bot_{j=0}^k s_j \G_j'
\end{equation}
where $\G_j' = \BZ_{p^{k+j}} \oplus \BZ_{p^{k-j}}$. By previous arguments, $\G'$ is the
group of fusion rules corresponding to a suitable pair $(G',\w')$. Thus (\ref{eq12.4}) and
(\ref{eq12.5}) enjoy the same relation to $(G',\w')$ that (\ref{eq12.2}) and (\ref{eq12.3})
do to $(G,\w)$.

Thus the canonical bijection (duality) which relates $\G$ and $\G'$ is achieved via the
duality between the partitions $\l$ and $\l^t$ which are associated to them. This explains
part (ii) of the theorem.
\begin{example}
{\rm
In Table 1 below, we illustrate this duality in case the two groups
in question are $G=\BZ_{p^2} \times \BZ_{p^2}\times
\BZ_{p^2}$ and $G'=\BZ_{p^3} \times \BZ_{p^3}$. Groups $\G$, $\G'$ on the same line are in
duality.
\begin{center}
\begin{tabular}{|c|c||c|c|c|}
\hline \multicolumn{2}{|c||} {$G=\BZ_{p^2} \times \BZ_{p^2}\times
\BZ_{p^2}$} & \multicolumn{2}{c|}{$G'=\BZ_{p^3}
\times\BZ_{p^3}$} & \\ \hline $(m_1,m_2)$ & $\G$ &
$(s_1,s_2,s_3)$ & $\G'$ & $f$\\ \hhline{|=|=#=|=|=|}
$(0,3)$ &
$\BZ_{p^4} \times \BZ_{p^4}\times \BZ_{p^4}$ & $(0,0,2)$ & $\BZ_{p^6}
\times \BZ_{p^6}$ & 1\\ \hline
$(1,2)$ & $\BZ_{p^4} \times \BZ_{p^4}\times \BZ_{p^3}\times \BZ_{p}$ & $(0,1,1)$ & $\BZ_{p^6}
\times \BZ_{p^5}\times \BZ_{p}$& 2\\ \hline
$(0,2)$ & $\BZ_{p^4} \times \BZ_{p^4}\times \BZ_{p^2}\times \BZ_{p^2}$ &
$(0,2,0)$ & $\BZ_{p^5} \times \BZ_{p^5}\times \BZ_{p}\times \BZ_{p}$& 1\\ \hline
$(2,1)$ & $\BZ_{p^4} \times \BZ_{p^3}\times \BZ_{p^3}\times \BZ_{p}\times \BZ_{p} $ &
$(1,0,1)$ &
$\BZ_{p^6}
\times \BZ_{p^4}\times \BZ_{p^2}$& 2\\ \hline
$(1,1)$ & $\BZ_{p^4} \times \BZ_{p^3}\times \BZ_{p^2}\times \BZ_{p^2}\times \BZ_{p} $ &
$(1,1,0)$ & $\BZ_{p^5} \times \BZ_{p^4}\times \BZ_{p^2}\times \BZ_{p}$& 2\\ \hline
$(0,1)$ & $\BZ_{p^4} \times \BZ_{p^2}\times \BZ_{p^2}\times \BZ_{p^2}\times \BZ_{p^2} $ &
$(0,2,0)$ & $\BZ_{p^4} \times \BZ_{p^4}\times \BZ_{p^2}\times
\BZ_{p^2} $& 1\\ \hline
$(2,0)$ & $\BZ_{p^3} \times \BZ_{p^3}\times \BZ_{p^2}\times \BZ_{p^2}\times \BZ_{p} \times
\BZ_{p}$ & $(0,1,0)$ & $\BZ_{p^5} \times \BZ_{p^3}\times \BZ_{p^3}\times \BZ_{p}$& 1\\ \hline
$(3,0)$ & $\BZ_{p^3} \times \BZ_{p^3}\times \BZ_{p^3}\times \BZ_{p}\times \BZ_{p} \times
\BZ_{p}$ & $(0,0,1)$ & $\BZ_{p^6} \times \BZ_{p^3}\times \BZ_{p^3}$& 1\\ \hline
$(1,0)$ & $\BZ_{p^3} \times \BZ_{p^2}\times \BZ_{p^2}\times \BZ_{p^2}\times \BZ_{p^2} \times
\BZ_{p}$ & $(1,0,0)$ & $\BZ_{p^4} \times \BZ_{p^3}\times \BZ_{p^3}\times \BZ_{p^2}$& 1\\ \hline
$(0,0)$ & $\left(\BZ_{p^2}\right)^6$ &
 $(0,0,0)$ & $\left(\BZ_{p^3}\right)^4$ & 0\\ \hline
\end{tabular}\\
Table: 1
\end{center}
}
\end{example}

Turning to part (iii), we first establish:
\begin{lem}\label{l12.1}
With the notation of (\ref{eq12.2}), let $f$ be the number of indices $a \ge 1$ such that
$m_a \ge 1$. Then the number of inequivalent, non-degenerate bilinear forms $b$ on $\G$
with metabolizer (isomorphic to) $G$ is $2^f$.
\end{lem}
\pf From Theorem \ref{t11.5}, we only need to enumerate the equivalence classes of
bilinear forms that arise as in (\ref{eq12.1}) and (\ref{eq12.2}) above.

To show that this number is $2^f$, let us assume for convenience that $p \equiv 1 \pmod{4}$
(the case $p \equiv 3 \pmod{4}$ is proved similarly). According to \cite{Wa63} and the
discussion following equation (\ref{gs:eq1}), the only forms $\G_a$ with cyclic metabolizer
are of type $A \oplus A$ or type $B \oplus B$ if $1 \le a \le n-1$; type $A$ or type $B$ if
$a=n$; type $A \oplus A$ ($\cong$ type $B \oplus B$) if $a=0$. Furthermore, if $a \ge 1$
the only forms on
$m_a \G_a$ with metabolizer $\left(\BZ_{p^n}\right)^{m_a}$ are (up to equivalence) of the shape
$m_a(A \oplus A)$ or $(m_a -1)(A \oplus A) \oplus (B \oplus B)$.Thus, as long as $m_a \ge
1$ and $a \ge 1$, there are in any case just two  inequivalent forms on $m_a\G_a$ with
homogeneous metabolizer
$\left(\BZ_{p^n}\right)^{m_a}$. Now the lemma follows immediately. \qed\\

Note that the integer $f$ of Lemma \ref{l12.1} is the number of unequal parts of the
partition
$\l$. But it is easy to see (and well-known) that this is equal also to the number of
unequal parts of $\l^t$. So the canonical bijection of part (ii) of Theorem \ref{t12.2}
also induces bijections between the inequivalent braided monoidal categories of the form
$\M{\H}$ and
$\M{D^{\w'}(G')}$ which have  fusion rules $\G$, $\G'$ respectively. This holds for all
such $\G$, so that the {\em number} of inequivalent braided tensor categories $\M{\H}$ is
the same as the numbers of type $\M{D^{\w'}(G')}$.

This is not quite a natural bijection. However, we saw above that the allowable bilinear
forms on $\G$ are naturally indexed by $f$-tuples consisting of $A$'s and $B$'s. More
precisely, let us again assume for convenience that $p\equiv 1\pmod{4}$. If $a \ge 1$ with
$m_a
\ge 1$ then we may take all but one $\G_a$ in (\ref{eq12.2}) to be type $A$ or
$A
\oplus A$ (according to whether $a=n$ or $a < n$); the remaining $\G_a$ may be taken to be {\em
either} type $A$ {\em or} type $B$ (or type $A \oplus A$ or $B \oplus B$). All such $f$-tuples
are allowable, so that there are exactly $2^f$ such tuples, as claimed.  
  Since (by the lemma) this analysis applies not only to $\G$ but to the dual group
$\G'$, we obtain a ``natural'' bijection between bilinear forms, and hence braided monoidal
categories, by associating those forms which correspond to the same $f$-tuple of $A$'s and
$B$'s. This completes our discussion of part (iii) of the theorem.

Finally, the formulae in example \ref{ex12.2} follow from our analysis without difficulty.
For example, take the case $k=2$, so that $G \cong \BZ_{p^n} \oplus \BZ_{p^n}$. The possible
$(n+1)$-tuples $(m_0,m_1,\dots,m_n)$ with $\sum_{a \ge 0}m_a =2$ are trivially enumerated,
and each give rise to $2^f$ inequivalent forms as in Lemma \ref{l12.1}. If $m_0=2$ then $f=0$;
if $m_0=1$ or $m_0 =0$ and $m_a =2$ for some $a \ge 1$ then $f=1$; otherwise $f=2$. Thus the
total number of inequivalent allowable forms (or braided monoidal categories) is
$1+2n+2n+4{n \choose 2}=2n^2+2n+1$.
\begin{example}
{\rm
We can read off further examples of symmetry from table 1. Consider the
group
$\G =
\BZ_{p^4}\times
\BZ_{p^4}\times
\BZ_{p^2}\times
\BZ_{p^2}$: we see that the two groups $G = \BZ_{p^2}\times \BZ_{p^2}\times
\BZ_{p^2}$ and $G' = \BZ_{p^3}\times \BZ_{p^3}$ each admit $\G$ as the corrresponding
fusion rules for some choice of quadratic form. Indeed, the previous discussion shows
that there are exactly two non-degenerate quadratic forms on $\G$ with
metabolizers equal to both $G$ and $G'$. By Theorem \ref{t8.4} it follows that the
following holds: there are cohomology classes $[\w_1]$, $[\w_2]$ in $H^3(G,\BC^*)_{ab}$ and
$[\w_1']$,
$[\w_2']$ in $H^3(G',\BC^*)_{ab}$ such that $D^{\w_i}(G)$ and $D^{\w'_i}(G')$, $i=1,2$, are
gauge equivalent.
}
\end{example}

\bibliographystyle{amsalpha}
%\bibliography{mybibl}
\providecommand{\bysame}{\leavevmode\hbox to3em{\hrulefill}\thinspace}

\end{document}